\documentclass[12pt]{article} 
\usepackage[sectionbib]{natbib}
\usepackage{array,epsfig,fancyheadings,rotating}
\usepackage[unicode=true,colorlinks,urlcolor=blue,linkcolor=blue,citecolor=blue,bookmarksnumbered=true]{hyperref}
\usepackage{mathrsfs}
\usepackage[ruled]{algorithm2e}
\usepackage{sectsty, secdot}
\sectionfont{\fontsize{12}{14pt plus.8pt minus .6pt}\selectfont}
\renewcommand{\theequation}{\thesection\arabic{equation}}
\subsectionfont{\fontsize{12}{14pt plus.8pt minus .6pt}\selectfont}

\usepackage{amsmath}
\allowdisplaybreaks
\usepackage{amssymb}
\usepackage{amsfonts}
\usepackage{multirow}
\usepackage{amsthm}
\usepackage{bm}
\usepackage{scalerel,stackengine}
\stackMath
\newcommand\widecheck[1]{%
\savestack{\tmpbox}{\stretchto{%
\scaleto{%
  \scalerel*[\widthof{\ensuremath{#1}}]{\kern-0pt\bigwedge\kern-0pt}%
  {\rule[-\textheight/2]{0.5ex}{\textheight}}
}{\textheight}%
}{0.5ex}}%
\stackon[1.2pt]{#1}{\scalebox{-0.8}{\tmpbox}}%
}

\newtheorem{thm}{Theorem}

\newcommand{\R}{\mathbb{R}}
\def\diag{{\mathrm {diag}}}
\def\w{{\mathrm {w}}}

\begin{document}


\renewcommand{\baselinestretch}{1.4}





\fontsize{12}{14pt plus.8pt minus .6pt}\selectfont \vspace{0.8pc}
\centerline{\large\bf Online Generalized Additive Model }
\vspace{.4cm} 
\centerline{Ying Yang and Fang Yao$^*$\footnote{Corresponding author: Fang, Yao (fyaomath@pku.edu.cn).} } 
\vspace{.4cm} 
\centerline{\it  Department of Probability and Statistics, School of Mathematical Sciences,}
\centerline{\it  Center for Statistical Science, Peking University, Beijing, China}
\vspace{.55cm} \fontsize{9}{11.5pt plus.8pt minus.6pt}\selectfont


\begin{quotation}
\noindent {\it Abstract:}
Additive models and generalized additive models are effective semiparametric tools for multidimensional data. 
In this article we propose an online smoothing backfitting method for generalized additive models  with local polynomial smoothers. 
The main idea is to use a second order expansion to approximate the nonlinear integral equations to  maximize the local quasilikelihood and store the coefficients as the sufficient statistics which can be updated in an online manner by a dynamic candidate bandwidth method.
The updating procedure only depends on the stored sufficient statistics and the current data block. 
We derive the asymptotic normality as well as the relative efficiency lower bounds of the online estimates, which provides insight into the relationship between estimation accuracy and computational cost driven by the length of candidate bandwidth sequence. Simulations and real data examples are provided to validate our findings.
\\

\vspace{9pt}
\noindent {\it Key words and phrases:}
Generalized additive models; Online learning; Streaming data; Efficiency
\par
\end{quotation}\par

\def\thefigure{\arabic{figure}}
\def\thetable{\arabic{table}}

\renewcommand{\theequation}{\thesection.\arabic{equation}}

\fontsize{11}{14pt plus.8pt minus .6pt}\selectfont

\section{Introduction}

Additive models (AM) and generalized additive models (GAM) are important tools in nonparametric regression for dimension reduction \citep{1989Linear,1990Generalized}. A number of methods have been developed to fit AM and GAM, such as the classical backfitting algorithm \citep{1989Linear,opsomer1997fitting,2000Asymptotic} and the smooth backfitting \citep{mammen_existence_1999} for AM, the local scoring backfitting \citep{1990Generalized}, the local likelihood \citep{kauermann2003local} and the smooth backfitting \citep{yu_smooth_2008} for GAM. These approaches all depend on the backfitting algorithm. Hence the computation would be costly when the sample size becomes large, which limits their applications in real world. 

In this work, we propose an online method for fitting AM and GAM which processes the data in a streaming fashion and update the model without storing the previous data.
Such online methods can be used to process the out-of-memory data and output the results in real-time, and have been prevalent in machine learning and statistics. 
Various stochastic gradient descent problems have been extended to the streaming settings, see \cite{Langford2009Sparse}, \cite{Duchi2009Efficient}, \cite{Xiao2010Dual} and \cite{Dekel2012Optimal}. \cite{lin2011aggregated} proposed an aggregated estimating equation for generalized linear models. \cite{Schifano2016Online} studied the online predictive residual test for linear models and improved the aforementioned method of \cite{lin2011aggregated}.  Classification based on linear and quadratic discriminant analyses have also been studied, see \cite{hiraoka2000convergence}, \cite{kim2007incremental} and \cite{pang2005incremental}. Following the custom of these papers, we refer to the classical approaches using the full data as {\em batch} methods. To the best of our knowledge, we are the first to study AM and GAM in the online context.

We focus on the smooth backfitting method for its well-developed theoretical properties. In specific, 
the classical backfitting requires that the covariates are nearly independent,
the local scoring estimate by \cite{1990Generalized} is defined implicitly and the theoretical properties are not well
understood,
and the local likelihood estimate by \cite{kauermann2003local} assumes that the algorithm is convergent.
In contrast, the smooth backfitting algorithm converges with theoretical guarantees and attains the oracle bias and variance when a local linear smoother is adopted with no assumption on the dependence between the covariates. 

Note that AM can be viewed as a special case of GAM, i.e., a GAM with the identity link function, we focus on the latter in this work. The main idea of the classical smooth backfitting for GAM is to maximize a smoothed quasilikelihood which is a natural extension of parametric quasilikelihood estimation. This maximization need solve a system of nonlinear integral equations by a double iteration scheme. Each outer iteration is a linear approximation of the smoothed quasilikelihood equations and is equivalent to a projection onto a Hilbert space equipped with a smoothed squared norm which is calculated by the inner loop. The main obstacle to conduct smooth backfitting onlinely is that the estimation equations are nonlinear and there is no ad hoc statistic. We propose to use a second order expansion and store the coefficients as the sufficient statistics. These statistics depend on the data-driven bandwidth which is changing during the data collection. We employ the dynamic candidate bandwidth method proposed by \cite{yang2021ofda} to approximate the optimal bandwidth by a sequence of dynamic candidates and combine the corresponding statistics across blocks to update the estimate. 

We derive the asymptotic normality of the proposed estimates. When compared to the batch competitors using the full data, the online estimates of the component functions have an explicit lower bound of the relative efficiency in terms of asymptotic integrated mean squared errors which is proportional to the length of candidate bandwidth sequence. This bound is the same as the one-dimensional online nonparametric regression. 
When the length is larger than 10, the efficiency is higher than 98\%.
For computational complexity, the classical batch algorithm need store and iterate $O(N_K)$ data, and our proposed method only need $O(n_K)$, where $K$ is the number of full blocks, $N_K$ is the full sample size up to block $K$ and $n_K$ is the subsample size of the $K$th block. In online problems, $n_K\ll N_K$, and hence the proposed online method attains desirable estimation accuracy with considerably lower computing cost.

The rest of the paper is organized as follows. In Section \ref{sec:method}, we first review the classical smooth backfitting and then present the online smooth backfitting method. We also introduce the dynamic candidate bandwidth method for kernel-based estimates and the online bandwidth selection. In Section \ref{sec:alg}, we delineate the computing details as well as a detailed algorithm. In Section \ref{sec:th}, we study the asymptotic properties of the online estimates and the relative efficiency compared to the classical batch ones, and prove the convergence of the algorithm. The numerical experiments including simulations and real data applications are displayed in Section \ref{sec:numerical}, which further verifies the usefulness of our method.

\section{Methodology}
\label{sec:method}
\subsection{Classical smooth backfitting}

Let $Y$ be a random variable and $\bm{X}=(X_1,\ldots,X_d)^\top$ be a random vector of dimension $d$. 
Assume that $(\bm{X}_1,Y_1),\ldots,(\bm{X}_N,Y_N)$ is a random sample drawn from $(\bm{X},Y)$.
In GAM, the conditional mean of the response given a $d$-dimensional covariate vector $m(\bm{x})=E(Y\mid \bm{X}=\bm{x})$ is modeled via a known link $g$ by a sum of unknown component functions,
\begin{equation}
\label{model gam}
g(m(\bm{x}))=\beta_0+\beta_1(x_1)+\cdots+\beta_d(x_d).
\end{equation}
By employing a suitable link $g$, it allows wider applicability than ordinary additive models.

We now briefly review the classical smoothed backfitting method. 
Let $\beta^*(\bm{x})=\beta_0^*+\sum_{j=1}^d\beta_j^*(x_j)$ be the true additive function and let $\beta_{1j}^*(x_j)=h_j\beta_j^{\prime}(x_j)$ for $1\leq j\leq d$. Further denote $\bm{\beta}=(\beta,\beta_{11},\ldots,\beta_{1d})$ with $\beta(\bm{x})=\beta_0+\sum_{j=1}^d\beta_j(x_j)$, and define
\begin{equation*}
\beta(\bm{u},\bm{x})=\beta(\bm{x})+\sum_{j=1}^d\left(\frac{u_j-x_j}{h_j}\right)\beta_{1j}(x_j).
\end{equation*}
Suppose that the conditional variance is modeled as
\begin{equation}\label{eq:V}
var(Y |\bm{X} = \bm{x}) = V (m(\bm{x}))
\end{equation}
for some positive function $V$. The smoothed quasilikelihood is an estimator of the expected quasilikelihood $E(Q(g^{-1}(\bm{\beta}), Y )\mid  \bm{X} = \bm{x})$ given by
\begin{equation*}
\widehat{SQ}(\bm{\beta})=\int N^{-1}\sum_{i=1}^{N}Q(g^{-1}(\beta(\bm{X}_i,\bm{x})),Y_i)\mathcal{K}_{\widehat{\bm{h}}}(\bm{x},\bm{X}_i)d\bm{x},
\end{equation*}
where $Q$ is the quasilikelihood function with $\partial Q(m,y)/\partial m=(y-m)/V(m)$, $\mathcal{K}_{\widehat{\bm{h}}}(\bm{x},\bm{X}_i)=\prod_{j=1}^d \mathcal{K}\big((X_{ij}-x_j)/\widehat{h}_{j}\big)$ is the kernel function and $\widehat{\bm{h}}=(\widehat{h}_{1},\ldots,\widehat{h}_{d})^\top$ is the bandwidth.

Suppose that $Q(g^{-1}(\bm{\beta}),y)$ is strictly concave as a function of $\bm{\beta}$ for each $y$. Then $Q(g^{-1}(\bm{\beta}),y)$ satisfies that the (conditional) Bartlett identities is not monotone in $\bm{\beta}$ for every $ \bm{x}$. Thus $\widehat{SQ}$ has a unique maximizer $\widehat{\bm{\beta}}$ with probability tending to one which satisfies
\begin{align}
\label{eq:goal}
d\widehat{SQ}(\bm{\beta};\bm{f})&=0\nonumber\\
\textrm{for all }\bm{f}=(f,f_{11},\ldots,f_{1d})\textrm{ with }f(\bm{x})&=f_0+\sum_{j=1}^df_j(x_j)\textrm{ and }f_{1j}(\bm{x})=f_{1j}(x_j),
\end{align}
where $dSQ(\bm{\beta};\bm{f})$ is the Fréchet differential of the functional $SQ$ at $\bm{\beta}$ with increment $\bm{f}$. As discussed in \cite{yu_smooth_2008}, $\widehat{\bm{\beta}}$ can be viewed as the projection of the full dimensional local linear estimator on to an appropriate Hilbert space.
Specifically, let 
$\bm{X}_{i}(\bm{x})=(1,X_{i1}-x_1,\ldots,X_{id}-x_d)^\top$, 
$\bm{X}(\bm{x})=(\bm{X}_{1}(\bm{x}),\ldots,\bm{X}_{N}(\bm{x}))^\top$,
$\bm{Y}=(Y_1,\cdots,Y_N)^\top$ and $\bm{\mathcal{K}}_N(\bm{x};\widehat{\bm{h}})=\diag\{\mathcal{K}_{\widehat{\bm{h}}}(\bm{x},\bm{X}_1),\cdots,\mathcal{K}_{\widehat{\bm{h}}}(\bm{x},\bm{X}_N)\}$. 
Let $\mathcal{F}_0=\{\bm{f}=(f,f_{11},\ldots,f_{1d}):\R^{d}\rightarrow\R^{d+1}\}$ and $\mathcal{F}_0(\bm{V})$ be the space $\mathcal{F}_0$ equipped with the norm 
\begin{equation}
\label{norm}
\Vert \bm{f}\Vert_{\bm{V}}=\left\{\int\bm{f}(\bm{x})^\top\bm{V}(\bm{x})\bm{f}(\bm{x})d\bm{x}\right\}^{1/2},
\end{equation}
where $\bm{V}(\bm{x})$ is a $(d+1)\times(d+1)$ matrix.
Define the Hilbert space
\begin{align}
\label{space H}
\mathcal{H}(\bm{V})=\Big\{&\bm{f}=(f,f_{11},\ldots,f_{1d})\in\mathcal{F}_0(\bm{V}): f(\bm{x})=f_0+\sum_{j=1}^df_j(x_j),\textrm{ where }\nonumber\\
 &f_j:\R\rightarrow\R, f_{1j}(\bm{x})=g_j(x_j),\textrm{ where }g_j:\R\rightarrow\R, j=1,\ldots,d\Big\},
\end{align}
and define
$\widehat{\bm{V}}_N(\bm{x},\widehat{\bm{\beta}};\widehat{\bm{h}})=N^{-1}\bm{X}(\bm{x})\diag[\nu_{1}(\bm{x},\widehat{\bm{\beta}};\widehat{\bm{h}}),\ldots,\nu_{N}(\bm{x},\widehat{\bm{\beta}};\widehat{\bm{h}})]\bm{X}(\bm{x})^\top$, where $\nu_{i}(\bm{x},\bm{\beta};\widehat{\bm{h}})=q_2(\beta(\bm{X}_{i},\bm{x}),Y_{ki})\mathcal{K}_{\widehat{\bm{h}}}(\bm{x},\bm{X}_{i})$.
Suppose that $\widehat{\bm{\beta}}_{full}$ is the full dimensional local linear solution to $d\widehat{SQ}(\bm{\beta};\bm{f})=0$. 
Then the solution of \eqref{eq:goal} can be viewed as the projection of full dimension local linear estimate $\widehat{\bm{\beta}}_{full}$ onto $\mathcal{H}(\widehat{\bm{V}}_N(\bm{x},\widehat{\bm{\beta}};\widehat{\bm{h}}))$.

Take the ordinary additive model as an instance, which corresponds to a GAM with the identity link, i.e., $g(x)=x$. The underlying model becomes
\begin{equation*}
m(\bm{x})=\beta(\bm{x})=\beta_0+\beta_1(x_1)+\cdots+\beta_d(x_d),
\end{equation*}
and the full dimensional local linear estimate is
\begin{equation*}
\widehat{m}_{full}(\bm{x})=e_{1}^\top\left\{\bm{X}(\bm{x})^\top\bm{\mathcal{K}}_N(\bm{x};\widehat{\bm{h}})\bm{X}(\bm{x})\right\}^{-1}\left\{\bm{X}(\bm{x})^\top\bm{\mathcal{K}}_N(\bm{x};\widehat{\bm{h}})\bm{Y}\right\},
\end{equation*}
where $\bm{\mathcal{K}}_N(\bm{x};\widehat{\bm{h}})=\diag\{\mathcal{K}_{\widehat{\bm{h}}}(\bm{x},\bm{X}_{1}),\cdots,\mathcal{K}_{\widehat{\bm{h}}}(\bm{x},\bm{X}_{N})\}$.
With $\nu_{i}(\bm{x},\bm{\beta};\widehat{\bm{h}})=\mathcal{K}_{\widehat{\bm{h}}}(\bm{x},\bm{X}_{i})$ and $\widehat{\bm{V}}_N(\bm{x},\widehat{\bm{\beta}};\widehat{\bm{h}})$ defined accordingly, the  smooth backfitting estimate for AM is the projection of $\widehat{m}_{full}(\bm{x})$ onto $\mathcal{H}(\widehat{\bm{V}}_N(\bm{x},\widehat{\bm{\beta}};\widehat{\bm{h}}))$, as first studied in \cite{mammen_existence_1999}.

\subsection{Online smooth backfitting}

In the online context, assume that we observe the $k$th data block at time $k$ which contains $\bm{Y}_k=(Y_{k1},\ldots,Y_{kn_k})^\top$ and $\bm{X}_k=(\bm{X}_{k1},\ldots,\bm{X}_{kn_k})^\top$, where $\bm{X}_{ki}=(X_{ki1},\ldots,X_{kid})^\top$, and $(\bm{X}_{ki},Y_{ki})$ is a random sample drawn from $(\bm{X},Y)$ for $i=1,\ldots,n_k$. Denote the current time as $K$ and the terminal time as $K_{max}$ which may tend to infinity.
If all previous data are available, the estimated smoothed quasilikelihood by batch method is
\begin{equation}\label{eq:SQ}
\widehat{SQ}(\bm{\beta})=\int N_K^{-1}\sum_{k=1}^K\sum_{i=1}^{n_k}Q(g^{-1}(\beta(\bm{X}_{ki},\bm{x})),Y_{ki})\mathcal{K}_{\widehat{\bm{h}}_K}(\bm{x},\bm{X}_{ki})d\bm{x},
\end{equation}
where $\widehat{\bm{h}}_K$ is the bandwidth selected based on $K$ blocks.
For conciseness, we introduce the following notations. Recall that $n_k$ is the sample size of the $k$th data block and $N_K=\sum_{k=1}^Kn_k$ is the full sample size up to time $K$.
Denote the weight of sample size of the $k$th block at time $K$ as 
\begin{equation}\label{eq:weight}
\w_{k|K}=n_k/N_K.
\end{equation}
For $\bm{h}=(h_1,\ldots,h_d)^\top$ and $j=1,\ldots,d$, define the following functions,
{\begin{align}
\label{eq:F hat}
&{G}_{k0}(\bm{\beta};{\bm{h}})=n_k^{-1}\int \sum_{i=1}^{n_k}q_1\big(\beta(\bm{X}_{ki},\bm{x}),Y_{ki}\big)\mathcal{K}_{{\bm{h}}}(\bm{x},\bm{X}_{ki})d\bm{x},\nonumber\\
&{G}_{kj}(\bm{\beta};{\bm{h}})=n_k^{-1}\int\sum_{i=1}^{n_k}q_1\big(\beta(\bm{X}_{ki},\bm{x}),Y_{ki}\big)\mathcal{K}_{{\bm{h}}}(\bm{x},\bm{X}_{ki})d\bm{x}_{-j},\nonumber\\
&{G}_{k1j}(\bm{\beta};{\bm{h}})=n_k^{-1}\int\sum_{i=1}^{n_k}q_1\big(\beta(\bm{X}_{ki},\bm{x}),Y_{ki}\big)\left(\frac{X_{kij}-x_j}{{h}_{j}}\right)\mathcal{K}_{{\bm{h}}}(\bm{x},\bm{X}_{ki})d\bm{x}_{-j},\nonumber\\
&{\bm{G}}_k(\bm{\beta};{\bm{h}})=\big({G}_{k0},{G}_{k1},\ldots,{G}_{kd},{G}_{k11},\ldots,{G}_{k1d}\big)^\top(\bm{\beta};{\bm{h}}),\nonumber\\
&\bm{F}_{K}(\bm{\beta};{\bm{h}})=\sum_{k=1}^K\w_{k|K}\bm{G}_{k}(\bm{\beta};{\bm{h}}),
\end{align}}
where $q_\ell({\beta},y)$ is the $\ell$th derivative of $Q(g^{-1}(\beta),y)$ with respect to $\beta$ and $\bm{x}_{-j}$ denotes the vector $\bm{x}$ with the $j$ th component $x_j$ being deleted. 
Define
\begin{equation}\label{eq:nu}
\nu_{ki}(\bm{x},\bm{\beta};{\bm{h}})=q_2(\beta(\bm{X}_{ki},\bm{x}),Y_{ki})\mathcal{K}_{{\bm{h}}}(\bm{x},\bm{X}_{ki}).
\end{equation}
and
\begin{equation*}
\bm{V}_K(\bm{x},\bm{\beta};\widehat{\bm{h}}_K)=N_K^{-1}\bm{X}(\bm{x})\diag[\nu_{11}(\bm{x},\bm{\beta};\widehat{\bm{h}}_K),\ldots,\nu_{KN_K}(\bm{x},\bm{\beta};\widehat{\bm{h}}_K)]\bm{X}(\bm{x})^\top,
\end{equation*}
where
$\bm{X}(\bm{x})=(\bm{X}_{11}(\bm{x}),\ldots,\bm{X}_{KN_K}(\bm{x}))^\top$,
$\bm{X}_{ki}(\bm{x})=(1,X_{ki1}-x_1,\ldots,X_{kid}-x_d)^\top$.
For uniqueness, we add a norm constrain, then with $\widehat{SQ}$ defined in \eqref{eq:SQ}, \eqref{eq:goal} is equivalent to 
\begin{align}
\label{eq:constrain}
&{\bm{F}}_K(\bm{\beta};\widehat{\bm{h}}_K)=0,\nonumber\\
{\rm subject\ to\ }&\langle\bm{\beta}_j, e_1 \rangle_{\bm{V}_K}=\int\bm{\beta}_j^\top(\bm{x})\bm{V}_K(\bm{x};\widehat{\bm{h}}_K) e_1 d\bm{x}= 0,
\end{align}
where
$\bm{\beta}_j(\bm{x})=(\beta_{j}(x_j),0,\ldots,0,\beta_{1j}(x_j),0,\ldots,0)^\top\in \R^{d+1}$ and
$ e_j\in \R^{d+1}$ is a vector whose $j$th element is 1 and the others are 0.

Solving \eqref{eq:constrain} in an online manner is hindered by two problems. First, 
$\bm{G}_k(\bm{\beta};\widehat{\bm{h}}_K)$
is nonlinear with respect to $\bm{\beta}$, and the precise statistics do not exist. We propose to use a two-order expansion of $\bm{G}_k(\bm{\beta};\widehat{\bm{h}}_K)$ at $\widetilde{\bm{\beta}}_k$ to make the approximation, where $\widetilde{\bm{\beta}}_k$ is the estimate obtained at time $k$. Then the coefficients of the expansion are the statistics of interest. 
The second obstacle is that $\bm{G}_k(\bm{\beta};\widehat{\bm{h}}_K)$ depends on the bandwidth $\widehat{\bm{h}}_K$, which varies as $K$ tends to infinity. Storing the statistics at all $\widehat{\bm{h}}_K$ is prohibitive. We adopt the dynamic candidate bandwidth method proposed by \cite{yang2021ofda}, which generates a sequence of candidate bandwidths for each block and dynamically selects one to estimate $\widehat{\bm{h}}_K$. See details in Section \ref{sec: DCM}.
Denote the selected candidate of the $k$th block at time $K$ as $\widetilde{\bm{\eta}}_{k|K}$ which is dynamically changing as $K$ increases. 
Then $\bm{G}_{k}(\bm{\beta};\widehat{\bm{h}}_K)$ is approximated by $\bm{G}_{k}(\bm{\beta};\widetilde{\bm{\eta}}_{k|K})$.  

We now present the two-order expansion of $\bm{G}_{k}(\bm{\beta};\widetilde{\bm{\eta}}_{k|K})$ \eqref{eq:F hat}.
For any function $G(\bm{\beta};\bm{h})$, define its two-order expansion at $\widetilde{\bm{\beta}}$ as
{\small\begin{equation*}
T_2G(\bm{\beta};\bm{h})\Big|_{\widetilde{\bm{\beta}}}=G(\widetilde{\bm{\beta}};\bm{h})
+\left\{\frac{\partial G(\widetilde{\bm{\beta}};\bm{h})}{\partial\widetilde{\bm{\beta}}}\right\}^\top\left(\bm{\beta}-\widetilde{\bm{\beta}}\right)
+\frac{1}{2}\left(\bm{\beta}-\widetilde{\bm{\beta}}\right)^\top\left\{\frac{\partial^2G(\widetilde{\bm{\beta}};\bm{h})}{\partial\widetilde{\bm{\beta}}^2}\right\}\left(\bm{\beta}-\widetilde{\bm{\beta}}\right).
\end{equation*} }
Recall that $\widetilde{\bm{\beta}}_k$ is the online estimate obtained at time $k$.
We estimate ${\bm{F}}_K(\bm{\beta};\widehat{\bm{h}}_K)$ by
{\begin{align}
\label{eq:F tilde}
&\widetilde{\bm{F}}_{K}\bm{\beta}
=\sum_{k=1}^{K-1}\w_{k|K}\cdot T_2\bm{G}_k(\bm{\beta};\widetilde{\bm{\eta}}_{k|K})\Big|_{\widetilde{\bm{\beta}}_k}
+\w_{K|K}\cdot \bm{G}_{K}(\bm{\beta},\widetilde{\bm{h}}_K),
\end{align}}
where $\bm{G}_K$ is defined as in \eqref{eq:F hat} and $\w_{k|K}$ is defined as in \eqref{eq:weight}. 
To appreciate \eqref{eq:F tilde}, note that $\bm{G}_{K}(\bm{\beta},\widetilde{\bm{h}}_K)$ is available at time $K$, and $\sum_{k=1}^{K-1}\w_{k|K}\cdot T_2\bm{G}_k(\bm{\beta};\widetilde{\bm{\eta}}_{k|K})\Big|_{\widetilde{\bm{\beta}}_k}$ is a quadratic function of $\bm{\beta}$ whose coefficients are the aggregation of fixed block-wise statistics. Then formula \eqref{eq:F tilde} can update in a streaming manner as long as the coefficients of $\sum_{k=1}^{K-1}\w_{k|K}\cdot T_2\bm{G}_k(\bm{\beta};\widetilde{\bm{\eta}}_{k|K})\Big|_{\widetilde{\bm{\beta}}_k}$ are stored. Specifically, the statistics to store are the elements of
\begin{eqnarray}\label{eq:suff stat}
&\Omega_{K}^{suff}=\left\{\sum_{k=1}^{K}\w_{k|K}\bm{T}_{k|K}:\bm{T}_{k|K}\in\Omega_{k|K}^{sub}\right\},
\end{eqnarray} 
where $\Omega_{k|K}^{sub}$ is the set of sub-sufficient statistics based on each block for $k=1,\ldots,K$ as follows,
{\small\begin{eqnarray}\label{eq:sub suff stat}
&\Omega_{k|K}^{sub}=\left\{\bm{G}_{K}(\widetilde{\bm{\beta}}_k;\widetilde{\bm{\eta}}_{k|K}),\ 
\frac{\partial\bm{G}_{K}(\widetilde{\bm{\beta}}_k;\widetilde{\bm{\eta}}_{k|K})}{\partial\widetilde{\bm{\beta}}_k},\ 
\left(\frac{\partial\bm{G}_{K}(\widetilde{\bm{\beta}}_k;\widetilde{\bm{\eta}}_{k|K})}{\partial\widetilde{\bm{\beta}}_k}\right)^\top\widetilde{\bm{\beta}}_k,\ \right.\nonumber\\
&\left.
\frac{\partial^2\bm{G}_{K}(\widetilde{\bm{\beta}}_k;\widetilde{\bm{\eta}}_{k|K})}{\partial(\widetilde{\bm{\beta}}_k)^2},\ 
\left(\frac{\partial^2\bm{G}_{K}(\widetilde{\bm{\beta}}_k;\widetilde{\bm{\eta}}_{k|K})}{\partial(\widetilde{\bm{\beta}}_k)^2}\right)^\top\widetilde{\bm{\beta}}_k,\
(\widetilde{\bm{\beta}}_k)^\top\frac{\partial^2\bm{G}_{K}(\widetilde{\bm{\beta}}_k;\widetilde{\bm{\eta}}_{k|K})}{\partial(\widetilde{\bm{\beta}}_k)^2}\widetilde{\bm{\beta}}_k\right\}.
\end{eqnarray} }
We mention that in \eqref{eq:suff stat},  $\{\bm{T}_{k|K}\}_{k=1}^K$ are of the same form except the index $k$.

Recall that $\widetilde{\bm{\beta}}_K$ is the estimate at time $K$, i.e., the solution to $\widetilde{\bm{F}}_K\bm\beta=0$. It is proved in Theorem \ref{THM:AN} of Section \ref{sec:th} that $\widetilde{\bm{\beta}}_K$ converges to $\bm{\beta}^*$ as $K\rightarrow\infty$.
Note that $\bm{\beta}$ can not be solved explicitly from $\widetilde{\bm{F}}_K\bm\beta=0$, we propose to use the linear approximation, i.e., the updating equation at time $K$ for computing the $m$th iteration estimate $\widetilde{\bm{\beta}}_K^{(m)}$ is given by
\begin{equation}
\label{eq:linearapprox}
0=\widetilde{\bm{F}}_K\bm{\beta}^{(m-1)}+\widetilde{\bm{F}}_K^{\prime}\big(\bm{\beta}^{(m-1)}\big)\big(\bm{\beta}-\bm{\beta}^{(m-1)}\big).
\end{equation}
The computing details and an implemented algorithm are presented in Section \ref{sec:alg}. 
It is proved in Theorem \ref{THM:OUTER CONVERGENCE} and \ref{THM:INNER CONVERGENCE} of Section \ref{sec:th} that $\widetilde{\bm{\beta}}_K^{(m)}$ converges to $\widetilde{\bm{\beta}}_K$ as $m\rightarrow\infty$. 

\subsection{Dynamic candidate method}\label{sec: DCM}

We mention that statistics in $\Omega_{k|K}^{sub}$ \eqref{eq:sub suff stat} and $\Omega_{K}^{suff}$ \eqref{eq:suff stat} depend on the pseudo-bandwidths $\{\widetilde{\bm{\eta}}_{k|K}\}_{k=1}^K$. For each $k$, $\widetilde{\bm{\eta}}_{k|K}$ is selected from the candidate bandwidth sequence $\{\bm{\eta}_{k\ell}\}_{\ell=1}^L$ where $\bm{\eta}_{k\ell}=(\eta_{k1\ell},\cdots,\eta_{kd\ell})^\top$. Now we introduce the selection rule. 
Let $\widetilde{\bm{h}}_{k}$ be the estimated optimal bandwidth based on the asymptotic behavior of the proposed online method at time $k$, which decreases with respect to $k$ and is given explicitly in \eqref{eq:online h} of Section \ref{sec:h}. It is proved in Theorem \ref{THM:BAND} of Section \ref{sec:th} that $\widetilde{\bm{h}}_{k}$ is quite close to $\widehat{\bm{h}}_{k}$.
We set $\eta_{kj1}=\widetilde{h}_{kj}$ and $\eta_{kj1}> \eta_{kj2}> \ldots> \eta_{kjL}$ for $j=1,\ldots,d$. The explicit expression of $\eta_{kj\ell}$ is given in \eqref{eq:eta} of Section \ref{sec:th}.

\begin{figure}[htbp]
\centering
\includegraphics[width = 5.5in]{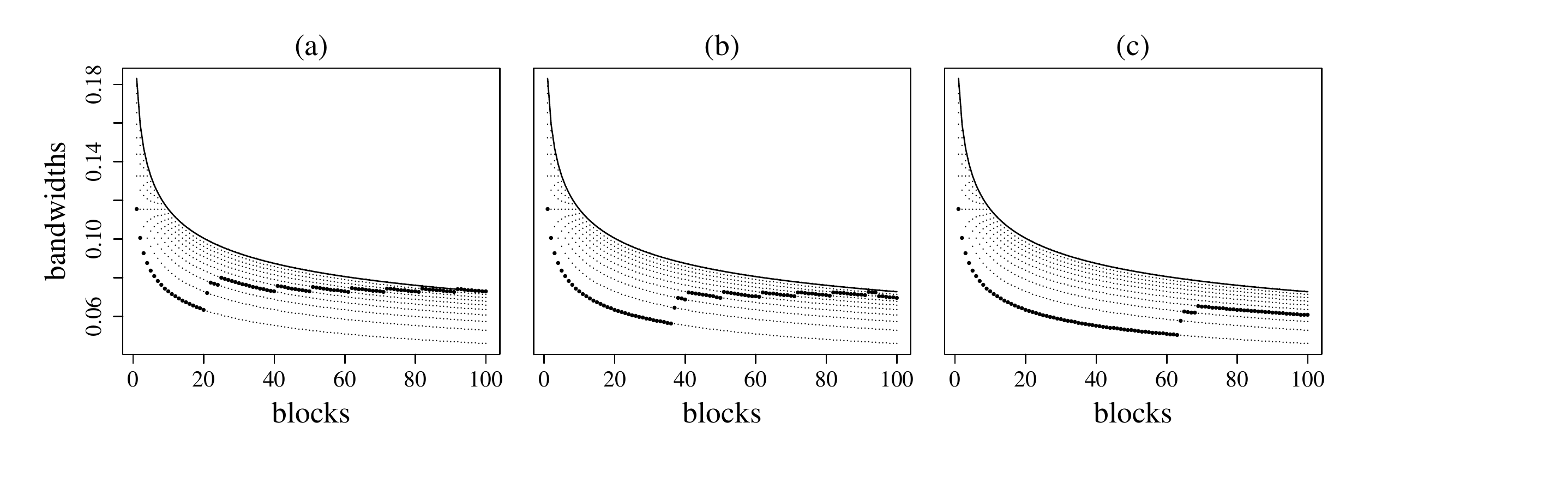}
\caption{\label{fig:sketch} A sketch for the dynamic candidate bandwidth method with $L=10$ and $K=100$ for the $j$th component, where the vertical dot series represent the dynamic candidate sequences $\{\eta_{kj\ell}\}_{\ell=1}^L$, $k=1,2,\ldots,K$, the bold points are the bandwidths used in ${\Omega}_{K\ell}^{suff}$ whose weighted average is the centroid $\phi_{Kj\ell}$ with $\ell=1,3,7$ corresponding to (a)-(c), and the solid line by connecting the largest values at each block corresponding to the online bandwidth estimates $\widetilde{h}_{kj},k=1,2,\ldots,K$.}
\end{figure}

We begin with computing the sets of sub-sufficient statistics with candidate bandwidths $\{\bm{\eta}_{K\ell }\}_{\ell=1}^L$ for the $K$th data block, denoted by $\{\Omega_{K\ell}^{sub}\}_{\ell=1}^L$, i.e.,
{\small\begin{eqnarray}\label{eq:sub suff seq}
&\Omega_{K\ell}^{sub}=\left\{\bm{G}_{K}(\widetilde{\bm{\beta}}_K;{\bm{\eta}}_{K\ell}),\ 
\frac{\partial\bm{G}_{K}(\widetilde{\bm{\beta}}_K;{\bm{\eta}}_{K\ell})}{\partial\widetilde{\bm{\beta}}_K},\ 
\left(\frac{\partial\bm{G}_{K}(\widetilde{\bm{\beta}}_K;{\bm{\eta}}_{K\ell})}{\partial\widetilde{\bm{\beta}}_K}\right)^\top\widetilde{\bm{\beta}}_K,\ \right.\nonumber\\
&\left.
\frac{\partial^2\bm{G}_{K}(\widetilde{\bm{\beta}}_K;{\bm{\eta}}_{K\ell})}{\partial(\widetilde{\bm{\beta}}_K)^2},\ 
\left(\frac{\partial^2\bm{G}_{K}(\widetilde{\bm{\beta}}_K;{\bm{\eta}}_{K\ell})}{\partial(\widetilde{\bm{\beta}}_K)^2}\right)^\top\widetilde{\bm{\beta}}_K,\
(\widetilde{\bm{\beta}}_K)^\top\frac{\partial^2\bm{G}_{K}(\widetilde{\bm{\beta}}_K;{\bm{\eta}}_{K\ell})}{\partial(\widetilde{\bm{\beta}}_K)^2}\widetilde{\bm{\beta}}_K\right\}.
\end{eqnarray} }   
We remark that $\Omega_{k|K}^{sub}$ \eqref{eq:sub suff stat} is dynamically selected from $\{\Omega_{k\ell}^{sub}\}_{\ell=1}^L$ according to $\widetilde{\bm{\eta}}_{k|K}$ for each $k$. 
For notation conciseness, let $\widetilde{\bm{T}}_{K\ell}$ be any element of $\Omega_{K\ell}^{sub}$. 
Then the combined statistic $\widetilde{\bm{T}}_{K\ell}^{C}$ is updated as follows,
\begin{eqnarray}
\label{eq:update T}
\widetilde{\bm{T}}_{K\ell}^{C}(x)&=(1-\w_{K|K})\widetilde{\bm{T}}_{K-1,j_{K\ell}}^{C}(x)+\w_{K|K}\widetilde{\bm{T}}_{K\ell}(x),
\end{eqnarray}
where $\w_{K|K}$ is the ratio of subsample size and the full sample size defined as in \eqref{eq:suff stat} and the index $j_{K\ell}$ is defined later.
To update $j_{K\ell}$ in an online fashion, we define the centroids $\{\bm\phi_{K\ell}\}_{\ell=1}^L$ as the weighted averages of all previous candidate bandwidths whose sub-sufficient statistics are aggregated, i.e.,
\begin{equation}
\label{eq:centroids}
\bm{\phi}_{K\ell}=(1-\w_{K|K})\bm{\phi}_{K-1,j_{K\ell}}+\w_{K|K}\bm{\eta}_{K\ell}.
\end{equation}
Then the index of $j_{K\ell}$ is defined as follows,
\begin{equation}\label{eq:d(l)}
j_{K\ell}={\rm argmin}_{i\in\{1,2,\ldots,L\}}\Vert\bm{\eta}_{K\ell}-\bm{\phi}_{K-1,i}\Vert_1,
\end{equation}
where the norm $\Vert \bm{u}-\bm{v}\Vert_1=\sum_{j=1}^d|u_j-v_j|$. Define
\begin{eqnarray}\label{eq:suff seq}
&\Omega_{K,\ell}^{suff}=\left\{\widetilde{\bm{T}}_{K\ell}^{C}:\widetilde{\bm{T}}_{K\ell}^{C}\textrm{ is updated as \eqref{eq:update T}}\textrm{ based on }\{\Omega_{K\ell}^{sub}\}_{K\ell}\right\}.
\end{eqnarray}   
With \eqref{eq:centroids} and \eqref{eq:d(l)}, the statistic $\widetilde{\bm{T}}_{K\ell}^{C}(x)$ employs the bandwidths closest to $\bm{\eta}_{K\ell}$ for all previous blocks (on average), as illustrated in Figure \ref{fig:sketch}, also see details in \cite{yang2021ofda}. 
Note that $\bm{\eta}_{K\ell}=\widetilde{\bm{h}}_K$, $\Omega_{K}^{suff}$ in \eqref{eq:suff stat} satisfies that $\Omega_{K}^{suff}=\Omega_{K1}^{suff}$. 
We emphasize that the algorithm only need store $L$ sets of statistics $\{\Omega_{K,\ell}^{suff}\}_{\ell=1}^L$ and $L$ candidate bandwidths $\{\bm{\eta}_{K\ell}\}_{\ell=1}^L$ throughout the procedure.

\subsection{Bandwidth selection}\label{sec:h}
We finally introduce the selection of $\widetilde{\bm{h}}_K=(\widetilde{h}_{K1},\cdots,\widetilde{h}_{Kd})^\top$. Denote
\begin{equation}\label{eq:sigma}
\sigma_j^2(x_j)=E\{V(m^*(\bm{X}))^{-1}g'(m^*(\bm{X}))^{-2}\mid X_j=x_j\}^{-1}.
\end{equation}
Let $p(\bm{x})$ be the density function of $\bm{X}$, and for $j=1,\ldots,d$, the marginal density is $p_j(x_j)=\int p(\bm{x})d\bm{x}_{-j}$.
From \cite{yu_smooth_2008}, one obtain the optimal bandwidth based on asymptotic integrated mean squared error at time $K$ is
\begin{equation}\label{eq:opt h}
h_{Kj}^*=C(\mathcal{K})\left\{\frac{\sigma_j^2}{\theta_j N_K}\right\}^{\frac{1}{5}},
\end{equation}
where $C(\mathcal{K})=\{R(\mathcal{K})/\mu_2(\mathcal{K})^2\}^{1/5}$ is a constant depends on the kernel function $K$ with $R(\mathcal{K})=\int \mathcal{K}(x)^2dx$ and $\mu_2(\mathcal{K})=\int x^2\mathcal{K}(x)dx$,
\begin{equation}\label{eq:theta sigma}
\sigma_j^2=\int\sigma_j^2(x)dx,\ \theta_j=\int \{(\beta_j^*){''}(x)\}^2p_j(x)dx.
\end{equation}
The online estimation of optimal bandwidth at time $K$ for the $j$th component is
\begin{equation}\label{eq:online h}
\widetilde{h}_{Kj}=C(\mathcal{K})\left\{\frac{\widetilde\sigma_{Kj}^2}{\widetilde\theta_{Kj} N_K}\right\}^{\frac{1}{5}},
\end{equation}
where $\widetilde{\theta}_{Kj}$ and $\widetilde{\sigma}_{Kj}^2$ are the online estimates of $\theta_j$ and $\sigma_j^2$, respectively. 

We suggest to use pilot estimates adopting the online method to approximate the unknown integral $\theta_j$ and the variance $\sigma_j^2$ for each $j=1,\ldots,d$.
Specifically, for $\sigma_j^2$, first fit another pilot online GAM with local linear smoother and candidate bandwidth sequence $\{\eta_{kjl}^{\sigma}\}_{\ell=1}^{L'}$, $k=1,2,\ldots,K$. Denote the corresponding estimate of $\beta_j^*$ by $\widetilde\beta_{Kj}^\sigma$ and estimate of $p_j$ by $\widetilde{p}_{Kj}=\int e_{1}^\top\{\bm{X}(\bm{x})^\top\bm{\mathcal{K}}_N(\bm{x};\widehat{\bm{h}})\bm{X}(\bm{x})\}e_{1}d\bm{x}_{-j}$. Then plug $\widetilde\beta_{Kj}^\sigma$ in \eqref{eq:sigma} to obtain $\widetilde{\sigma}_{Kj}^2$. 
For $\theta_j$, fit an online GAM with local quadratic smoother and candidate bandwidth sequence $\{\eta_{kjl}^{\theta}\}_{\ell=1}^{L'}$ to obtain the estimate $\widetilde\beta_{Kj}^{''}$, and integrate $\widetilde\beta_{Kj}^{''}\widetilde{p}_{Kj}$ to obtain $\widetilde{\theta}_{Kj}$.
We now present the selection of $\{\eta_{kjl}^{\sigma}\}_{\ell=1}^{L'}$ and $\{\eta_{kjl}^{\theta}\}_{\ell=1}^{L'}$. 
As discussed in \cite{yu_smooth_2008}, for local polynomials of different orders, the optimal bandwidths of GAM solved by smooth backfitting are of the same order as the one-dimensional case. 
The optimal bandwidths for estimating $\theta_j$ and $\sigma_j^2$ are $h_{Kj}^{\theta*}=G_j^*N_K^{-1/7}$ and $h_{Kj}^{\sigma*}=R_j^*N_K^{-2/9}$, respectively, where $R_j^*$ and $G_j^*$ involve unknown quantities depending on $\bm{\beta}^*$. 
We derive in Theorem \ref{THM:BAND} of Section \ref{sec:th} that when the estimated bandwidths satisfy
\begin{equation}
\label{eq:gr}
h_{Kj}^\theta=GN_K^{-1/7},\ h_{Kj}^\sigma=RN_K^{-1/5},\ G=O(1),\ R=O(1), 
\end{equation}
and the corresponding candidates satisfy
\begin{eqnarray}\label{eq:candidate r}
&\eta_{Kj1}^{\theta}=h_{Kj}^\theta,\ \eta_{Kj1}^{\theta}>\cdots>\eta_{KjL}^{\theta},\nonumber\\
&\eta_{Kj1}^{\sigma}=h_{Kj}^\sigma,\ \eta_{Kj1}^{\sigma}>\cdots>\eta_{KjL}^{\sigma},
\end{eqnarray}
the estimate $\widetilde{h}_{Kj}$ attains the same convergence rate $N_K^{2/7}$ as its batch competitor.
Based on extensive numerical experiments, we recommend to set $G$ and $R$ between $0.5$ and 1.

\section{Computational Details}\label{sec:alg}

To get an explicit form of the updating equation \eqref{eq:linearapprox}, define the following weights besides $\nu_{ki}$ in \eqref{eq:nu},
\begin{eqnarray}\label{eq:uw}
{u}_{ki}(\bm{x},\bm{\beta};\bm{\eta}_{k\ell})=q_1(\beta(\bm{X}_{ki},\bm{x}),Y_{ki};\bm{\eta}_{k\ell})\mathcal{K}_{\bm{\eta}_{k\ell}}(\bm{x},\bm{X}_{ki}),\nonumber\\
{\omega}_{ki}(\bm{x},\bm{\beta};\bm{\eta}_{k\ell})=q_3(\beta(\bm{X}_{ki},\bm{x}),Y_{ki};\bm{\eta}_{k\ell})\mathcal{K}_{\bm{\eta}_{k\ell}}(\bm{x},\bm{X}_{ki}).
\end{eqnarray}
Then let $\bm{X}_k(\bm{x})=(\bm{X}_{k1}(\bm{x}),\ldots,\bm{X}_{kn_k}(\bm{x}))^\top$, where $\bm{X}_{ki}\big(\bm{x})=\big(1,X_{ki1}-x_1,\ldots,X_{kid}-x_d\big)^\top$ for $i=1,\ldots,n_k$, and 
\begin{eqnarray*}
 \bm{D}_{u,k}(\bm{x},\bm{\beta};\bm{\eta}_{k\ell})=\diag[{u}_{k1}(\bm{x},\bm{\beta};\bm{\eta}_{k\ell}),\ldots,{u}_{kn_k}(\bm{x},\bm{\beta};\bm{\eta}_{k\ell})],\\
 \bm{D}_{\nu,k}(\bm{x},\bm{\beta};\bm{\eta}_{k\ell})=\diag[{\nu}_{k1}(\bm{x},\bm{\beta};\bm{\eta}_{k\ell}),\ldots,{\nu}_{kn_k}(\bm{x},\bm{\beta};\bm{\eta}_{k\ell})],\\
 \bm{D}_{\omega,k}(\bm{x},\bm{\beta};\bm{\eta}_{k\ell})=\diag[{\omega}_{k1}(\bm{x},\bm{\beta};\bm{\eta}_{k\ell}),\ldots,{\omega}_{kn_k}(\bm{x},\bm{\beta};\bm{\eta}_{k\ell})].
\end{eqnarray*}
With the above notations, define the following functions to generate the sub-sufficient statistics,
\begin{eqnarray}\label{eq:UVW hat}
&{\bm{U}}_k(\bm{x},\bm{\beta};\bm{\eta}_{k\ell})=n_k^{-1}\bm{X}_k(\bm{x})^\top \bm{D}_{u,k}(\bm{x},\bm{\beta};\bm{\eta}_{k\ell})\bm{X}_k(\bm{x}),\nonumber\\
&{\bm{V}}_k(\bm{x},\bm{\beta};\bm{\eta}_{k\ell})=n_k^{-1}\bm{X}_k(\bm{x})^\top \bm{D}_{\nu,k}(\bm{x},\bm{\beta};\bm{\eta}_{k\ell})\bm{X}_k(\bm{x}), \nonumber\\
&{\bm{W}}_k(\bm{x},\bm{\beta};\bm{\eta}_{k\ell})=n_k^{-1}\bm{X}_k(\bm{x})^\top \bm{D}_{\omega,k}(\bm{x},\bm{\beta};\bm{\eta}_{k\ell})\bm{X}_k(\bm{x}),\\
&{\bm{S}}_{kj}(\bm{x},\bm{\beta};\bm{\eta}_{k\ell})= {\bm{W}}_k(\bm{x},\bm{\beta};\bm{\eta}_{k\ell})\diag[(X_{k1j}-x_j),\ldots,(X_{kn_kj}-x_j)],\nonumber
\end{eqnarray}
where $j=1,\ldots,d$.
Then the sub-sufficient statistics evaluated at $\widetilde{\bm{\eta}}_{k\ell}$ for the $k$th block are 
\begin{eqnarray*}
&\widetilde{\bm{U}}_{k\ell}={\bm{U}}_k(\bm{x},\widetilde{\bm{\beta}}_k;\bm{\eta}_{k\ell}),\ 
\widetilde{\bm{V}}_{k\ell}={\bm{V}}_k(\bm{x},\widetilde{\bm{\beta}}_k;\bm{\eta}_{k\ell}),\\ 
&\widetilde{\bm{W}}_{k\ell}={\bm{W}}_k(\bm{x},\widetilde{\bm{\beta}}_k;\bm{\eta}_{k\ell}),\ 
\widetilde{\bm{S}}_{kj,\ell}={\bm{S}}_{kj}(\bm{x},\widetilde{\bm{\beta}}_k;\bm{\eta}_{k\ell}),\\
&\widetilde{\bm{Vb}}_{k\ell}=\bm{V}_{k\ell}\widetilde{\bm{\beta}}_k,\ 
\widetilde{\bm{Wb}}_{k\ell}=\bm{W}_{k\ell}\widetilde{\bm{\beta}}_k,\ 
\widetilde{\bm{bWb}}_{k\ell}=(\widetilde{\bm{\beta}}_k)^\top\bm{W}_{k\ell}\widetilde{\bm{\beta}}_k,\\ 
&\widetilde{\bm{Sb}}_{kj,\ell}=\bm{S}_{k\ell}\widetilde{\bm{\beta}}_k,\ 
\widetilde{\bm{bSb}}_{kj,\ell}=(\widetilde{\bm{\beta}}_k)^\top\bm{S}_{k\ell}\widetilde{\bm{\beta}}_k,
\end{eqnarray*}
i.e., we have the equivalent expression for $\Omega_{k\ell}^{sub}$ \eqref{eq:sub suff seq} for $\ell=1,\ldots,L$, 
\begin{align}\label{eq:sub seq2}
\Omega_{k\ell}^{sub}=\big\{&\widetilde{\bm{U}}_{k\ell},\widetilde{\bm{V}}_{k\ell},\widetilde{\bm{W}}_{k\ell},\widetilde{\bm{Vb}}_{k\ell},\widetilde{\bm{Wb}}_{k\ell},\widetilde{\bm{bWb}}_{k\ell},\widetilde{\bm{Sb}}_{kj,\ell},\widetilde{\bm{bSb}}_{kj,\ell}:j=1,\ldots,d\big\}.
\end{align}
For any element of $\Omega_{k\ell}^{sub}$, 
the corresponding aggregated version is defined the same as in \eqref{eq:update T}, and $\Omega_{K-1,\ell}^{suff}$ \eqref{eq:suff seq} can be expressed as follows,
\begin{align}
  \Omega_{K-1,\ell}^{suff}=\big\{&\widetilde{\bm{U}}_{K-1,\ell},\widetilde{\bm{V}}_{K-1,\ell},\widetilde{\bm{Vb}}_{K-1,\ell},\widetilde{\bm{W}}_{K-1,\ell},\widetilde{\bm{Wb}}_{K-1,\ell},\widetilde{\bm{bWb}}_{K-1,\ell},\nonumber\\
&\widetilde{\bm{S}}_{K-1,j\ell},\widetilde{\bm{Sb}}_{K-1,j\ell},\widetilde{\bm{bSb}}_{K-1,j\ell}:j=1,\ldots,d\big\}.\label{eq:suff seq2}
\end{align}
Note that at time $K$, ${\bm{U}}_K(\bm{x},\bm{\beta};\widetilde{\bm{h}}_{K})$ is aggregated with $\widetilde{\bm{U}}_{K-1,j_{K1}}$ to calculate $\widetilde{\bm{\beta}}_K$, where $j_{K\ell}$ are defined as in \eqref{eq:d(l)}.
For conciseness, we omit the subscript ``$K-1,j_{K1}$'' in the following description of the update procedure (for example, $\widetilde{\bm{U}}=\widetilde{\bm{U}}_{K-1,j_{K1}}$).
Further define
\begin{eqnarray*}
&\bm{D}_1^V(\bm{\beta})=\widetilde{\bm{V}}\bm{\beta}-\widetilde{\bm{Vb}},\ \bm{D}_1^W(\bm{\beta})=\widetilde{\bm{W}}\bm{\beta}-\widetilde{\bm{Wb}},\ \bm{D}_1^{S_j}=\widetilde{\bm{S}}_j\bm{\beta}-\widetilde{\bm{Sb}}_j,\\ 
&\bm{D}_2^W(\bm{\beta})=\bm{\beta}^\top\widetilde{\bm{W}}\bm{\beta}-2\bm{\beta}^\top\widetilde{\bm{Wb}}+\widetilde{\bm{bWb}},
\bm{D}_2^{S_j}(\bm{\beta})=\bm{\beta}^\top\widetilde{\bm{S}}_j\bm{\beta}-2\bm{\beta}^\top\widetilde{\bm{Sb}}_j+\widetilde{\bm{bSb}}_j,
\end{eqnarray*}
which are available functions of $\bm{\beta}$ with $\Omega_{K-1,j_{K1}}^{suff}$ stored.
Let $\widetilde{V}_{pq}$ be the $(p+1,q+1)$th element of the matrix $\widetilde{\bm{V}}=\widetilde{\bm{V}}_{K-1,j_{K1}}$, and $\widehat{V}_{K,pq}^{(m-1)}$ be the $(p+1,q+1)$th element of the matrix $\widehat{\bm{V}}_{K}(\bm{x},\bm{\beta}^{(m-1)};\widetilde{\bm{h}}_K)$. Define the following matrices of the $K$th block
\begin{align*}
&\widehat{\bm{M}}_{K,j}^{(m-1)}=
\begin{bmatrix}
\widehat{V}_{K,00}^{(m-1)} & 
\widehat{V}_{K,0j}^{(m-1)}\\
\widehat{V}_{K,0j}^{(m-1)} & 
\widehat{V}_{K,jj}^{(m-1)}\\
\end{bmatrix},
\quad
\widehat{\bm{M}}_{K,jl}^{(m-1)}=
\begin{bmatrix}
\widehat{V}_{K,00}^{(m-1)} & 
\widehat{V}_{K,0l}^{(m-1)}\\
\widehat{V}_{K,0j}^{(m-1)} & 
\widehat{V}_{K,jl}^{(m-1)}
\end{bmatrix},
\end{align*}
and the approximated version of previous blocks can be expressed as
\begin{align*}
&\widetilde{\bm{M}}_{K-1,j}^{(m-1)}=
\begin{bmatrix}
\widetilde{V}_{00}+e_{1}^\top \bm{D}_1^W(\bm{\beta}^{(m-1)})&
\widetilde{V}_{0j}+e_{j+1}^\top \bm{D}_1^W(\bm{\beta}^{(m-1)})\\
\widetilde{V}_{0j}+e_{1}^\top \bm{D}_1^{S_j}(\bm{\beta}^{(m-1)}) & 
\widetilde{V}_{jj}+e_{j+1}^\top \bm{D}_1^{S_j}(\bm{\beta}^{(m-1)})\\
\end{bmatrix},
\\
&\widetilde{\bm{M}}_{K-1,jl}^{(m-1)}=
\begin{bmatrix}
\widetilde{V}_{00}+e_{1}^\top \bm{D}_1^W(\bm{\beta}^{(m-1)}) &
\widetilde{V}_{0l}+e_{\ell +1}^\top \bm{D}_1^W(\bm{\beta}^{(m-1)})\\
\widetilde{V}_{0j}+e_{1}^\top \bm{D}_1^{S_j}(\bm{\beta}^{(m-1)})& 
\widetilde{V}_{jl}+e_{\ell +1}^\top \bm{D}_1^{S_j}(\bm{\beta}^{(m-1)})
\end{bmatrix},
\end{align*}
which depend only on the current value $\bm{\beta}^{(m-1)}$ and the stored statistics.
Further define the aggregations
\begin{align*}
&\widecheck{\bm{M}}_{K,j}^{(m-1)}=\int\left\{(1-\w_{K|K})\widetilde{\bm{M}}_{K-1,j}^{(m-1)}+\w_{K|K}\widehat{\bm{M}}_{K,j}^{(m-1)}\right\}d\bm{x}_{-j},\\
&\widecheck{\bm{M}}_{K,jl}^{(m-1)}=\int\left\{(1-\w_{K|K})\widetilde{\bm{M}}_{K-1,jl}^{(m-1)}+\w_{K|K}\widehat{\bm{M}}_{K,jl}^{(m-1)}\right\}d\bm{x}_{-(j,l)}.
\end{align*}
Then the outer iteration is
\begin{align}
\label{outer iter}
&\xi_{0}^{(m)}=-\left[\int\left\{(1-\w_{K|K})\left(\widetilde{V}_{00}+e_1^\top \bm{D}_1^W(\bm{\beta}^{(m-1)})\right)\xi +\w_{K|K}\widehat{V}_{K,00}^{(m-1)}\right\} d\bm{x}\right]^{-1}\nonumber\\
&\hspace*{1.5cm}\times\int\left\{(1-\w_{K|K})\left(\widetilde{U}_{00}+ e_1^\top \bm{D}_1^V(\bm{\beta}^{(m-1)})+ \bm{D}_2^W(\bm{\beta}^{(m-1)})\right)
+\w_{K|K}\widehat{U}_{K,00}^{(m-1)} \right\}d\bm{x},\nonumber\\
&\widecheck{\bm{M}}_{K,j}^{(m-1)}
\begin{bmatrix}
\xi_{j}^{(m)}(x_j)\\
\xi_{1j}^{(m)}(x_j)\\
\end{bmatrix}
=
\begin{bmatrix}
\widetilde\zeta_{j}^{(m)}(x_j)\\
\widetilde\zeta_{1j}^{(m)}(x_j)\\
\end{bmatrix}-\xi_{0}^{(m)}
\widecheck{\bm{M}}_{K,j}^{(m-1)}e_1-\int\sum_{\ell \neq j}\widecheck{\bm{M}}_{K,jl}^{(m-1)}
\begin{bmatrix}
\xi_{\ell }^{(m)}(x_\ell)\\
\xi_{1l}^{(m)}(x_\ell)\\
\end{bmatrix}
d\bm{x}_\ell,
\end{align}  
where
{\small\begin{align}\label{eq:tilde zeta}
&\widetilde\zeta_{j}^{(m)}(x_j)=-\int\left\{(1-\w_{K|K})\left(\widetilde{U}_{00}+ e_{1}^\top \bm{D}_1^V(\bm{\beta}^{(m-1)})+ \bm{D}_2^W(\bm{\beta}^{(m-1)})\right)
+\w_{K|K}\widehat{U}_{K,00}^{(m-1)} \right\}d\bm{x}_{-j},
\nonumber\\
&\widetilde\zeta_{1j}^{(m)}(x_j)=-\int\left\{(1-\w_{K|K})\left(\widetilde{U}_{0j}+e_{j+1}^\top \bm{D}_1^V(\bm{\beta}^{(m-1)})+ \bm{D}_2^{S_j}(\bm{\beta}^{(m-1)})\right)
+\w_{K|K}\widehat{U}_{K,0j}^{(m-1)} \right\}d\bm{x}_{-j}.
\end{align}}
Let $\widetilde{U}_{pq}$ be the $(p+1,q+1)$th elements of the matrix $\widetilde{\bm{U}}$. 
For the norm constrain, define
\begin{equation}\label{eq:c}
c_j^{(m)}=\left\{\int\left((1-\w_{K|K})\widetilde{V}_{00}^{(m-1)}+\w_{K|K}\widehat{V}_{K,00}^{(m-1)}\right)d\bm{x}\right\}^{-1}\times\left\{\int\widehat{A}_{j}^{(m-1)}dx_j\right\},
\end{equation}
where
\begin{align*}
\widehat{A}_{j}^{(m-1)}(x_j)=
&\int\left\{(1-\w_{K|K})\widetilde{V}_{00}+\w_{K|K}\widehat{V}_{K,00}^{(m-1)}\right\}d\bm{x}_{-j}\times\xi_j^{(m-1)}(x_j)\\
&+\int\left\{(1-\w_{K|K})\widetilde{V}_{0j}+\w_{K|K}\widehat{V}_{K,0j}^{(m-1)}\right\}d\bm{x}_{-j}\times\xi_{1j}^{(m-1)}(x_j).
\end{align*}
Then, it can be shown that the updating equation \eqref{eq:linearapprox} is equivalent to    
\begin{align}\label{eq:update beta}
&\widetilde{\beta}_0^{(m)}=\widetilde{\beta}_0^{(m-1)}+\widetilde{\xi}_0^{(m)}+\sum_{j=1}^dc_j^{(m)},\nonumber\\
&\widetilde{\beta}_j^{(m)}=\widetilde{\beta}_j^{(m-1)}+\widetilde{\xi}_j^{(m)}-c_j^{(m)},j=1,\ldots,d,\\
&\widetilde{\beta}_{1j}^{(m)}=\widetilde{\beta}_{1j}^{(m-1)}+\widetilde{\xi}_{1j}^{(m)},j=1,\ldots,d.\nonumber
\end{align}

Equation \eqref{outer iter} is solved by the following inner iteration, and the updating for the $j$ th component of the $r$th iteration cycle is given by
{\small\begin{align}
\label{inner iter}
\widecheck{\bm{M}}_{K,j}^{(m-1)}
\begin{bmatrix}
\xi_{j}^{[m,r]}(x_j)\\
\xi_{1j}^{[m,r]}(x_j)\\
\end{bmatrix}
=&
\begin{bmatrix}
\widetilde\zeta_{j}^{(m)}(x_j)\\
\widetilde\zeta_{1j}^{(m)}(x_j)\\
\end{bmatrix}-\xi_{0}^{(m)}
\widecheck{\bm{M}}_{K,j}^{(m-1)}e_1
-\int\sum_{\ell < j}\widecheck{\bm{M}}_{K,jl}^{(m-1)}
\begin{bmatrix}
\xi_{\ell }^{[m,r]}(x_\ell)\\
\xi_{1l}^{[m,r]}(x_\ell)\\
\end{bmatrix}
d\bm{x}_\ell\nonumber\\
&\hspace{1.8in}-\int\sum_{\ell > j}\widecheck{\bm{M}}_{K,jl}^{(m-1)}
\begin{bmatrix}
\xi_{\ell }^{[m,r-1]}(x_\ell)\\
\xi_{1l}^{[m,r-1]}(x_\ell)\\
\end{bmatrix}
d\bm{x}_\ell,
\end{align}  }
An implementation is given in Algorithm \ref{alg:ogam}. It is proved in Theorem \ref{THM:OUTER CONVERGENCE} and \ref{THM:INNER CONVERGENCE} of Section \ref{sec:th} that this double cycle algorithm converges to the solution of $\widetilde{\bm{F}}_K\bm{\beta}=0$ under regular conditions.

{\small\begin{algorithm}
\caption{Online GAM algorithm}\label{alg:ogam}
\LinesNumbered 
\While{$K\leq K_{max}$}{
  $\bm{\beta}^{(0)}=\widetilde{\bm{\beta}}_{K-1}$\;
  compute $\widetilde{\bm{h}}_{K}$ and $\bm{\eta}_{K\ell}$ as in \eqref{eq:online h} and \eqref{eq:eta}, respectively\;
  compute $j_{K\ell}$ as in \eqref{eq:d(l)}\;
  $m\leftarrow1$\;
  \While{not converge}
  {compute $\widehat{\bm{U}}_k(\bm{x}$, $\bm{\beta}^{(m-1)};\widetilde{h}_K)$, $\widehat{\bm{V}}_k(\bm{x},\bm{\beta}^{(m-1)};\widetilde{h}_K)$ and $\widehat{\bm{W}}_k(\bm{x},\bm{\beta}^{(m-1)};\widetilde{h}_K)$ as in \eqref{eq:UVW hat}\;
  $r\leftarrow1$\;
  \While{not converge}
  {update \eqref{inner iter} based on statistics in $\Omega_{K-1,j_{K1}}^{suff}$\;
  $r\leftarrow r+1$\;
  }
  compute $c_j^{(m)}$ and $\bm{\beta}^{(m)}$ as in \eqref{eq:c} and \eqref{eq:update beta}, respectively\;
  $m\leftarrow m+1$\;
  }
  $\widetilde{\bm{\beta}}_K\leftarrow\bm{\beta}^{(m)}$\;
  $K\leftarrow K+1$\;
  update the stored set of statistics $\big\{\Omega_{K-1,\ell}^{suff}\big\}_{\ell=1}^L$ as in \eqref{eq:suff seq2}\;
}
\end{algorithm}}

\section{Theoretical Analysis}
\label{sec:th}

We impose the following assumptions to conduct theoretical analysis.
\renewcommand{\theenumi}{A\arabic{enumi}}
\renewcommand{\labelenumi}{(\theenumi)}

\begin{enumerate}
\item\label{assump: density} The density of $\bm{X}$, i.e., $p(\bm{X})$, is bounded away from zero and infinity on its support $[0,1]^d$ and has continuous derivatives.
\item\label{assump: model} $q_2(u,y)<0$ for $u\in\R$ and $y$ in the range of the response, the link $g$ is strictly monotone and is three times continuously differentiable, $V$ defined in \eqref{eq:V} is strictly positive and twice continuously differentiable, and $\nu(\bm{x})=var(Y\mid \bm{X}=\bm{x})$ is continuous. The moment $E\vert Y\vert^{r_0}<\infty$ for some $r_0>5/2$.
\item\label{assump: components} The true component functions are twice continuously differentiable.
\item\label{assump: kernel} The kernel function $\mathcal{K}(\cdot)$ is a symmetric density function with compact support $[0,1]$.
\item\label{assump: pseudo-bandwidths} The bandwidths $\widetilde{h}_{Kj}=O_p(N_K^{-1/5})$ for $j=1,\ldots,d$
\item\label{assump:block} The block size satisfies $n_k/N_K\rightarrow0$, $k=1,2,\ldots,K$ as $K\rightarrow\infty$.
\end{enumerate}

Assumption \eqref{assump: density}-\eqref{assump: model} are general assumptions for generalized additive models. 
Assumption \eqref{assump: components}-\eqref{assump: kernel} are standard for local linear smoothing. 
Assumption \eqref{assump: pseudo-bandwidths} requires that the estimated bandwidths take the specific order $N_K^{-1/5}$. This is necessary to verify and simplify some complicated conditions as discussed in \cite{mammen_existence_1999}, and has become a common assumption for smooth backfitting algorithm, see \cite{mammen_bandwidth_2005} and \cite{yu_smooth_2008}.
Assumption \eqref{assump:block} indicates that the block size $n_k$ is relatively small compared to the full sample size $N_K$, which is natural in the online context.
To give a concise description, we define for $l\in\mathbb{Z}$,
\begin{equation}
\label{eq:eta moments}
\bm{\rho}_{Kl}=(\rho_{K1,l},\cdots,\rho_{Kd,l})^\top=\frac{1}{N_K}\sum_{k=1}^Kn_k\left(\widetilde{\bm{\eta}}_{k|K}\right)^l.
\end{equation}
We now state the asymptotic behavior of $\widetilde\beta_{Kj}(x_j)$, $j=1,\ldots,d$.

\begin{thm}
\label{THM:AN}
Under Assumption \eqref{assump: density}-\eqref{assump: pseudo-bandwidths}, recall that $\beta^*$ is the true function, for any $x_1,\ldots,x_d\in(0,1)$, then as $K\rightarrow\infty$, we have
\begin{align*}
diag\{v_j(x_j)\}^{-1/2}\big(\widetilde{\beta}_{K1}(x_1)-\beta_1^*(x_1)-b_1(x_1),\ldots,\widetilde{\beta}_{Kd}(x_d)-\beta_d^*(x_d)-b_d(x_d)\big)^\top\\
\overset{d}{\longrightarrow}N(0,I_d),&
\end{align*}
where $N(0,I_d)$ is the $d$-variate standard Normal distribution,
\begin{align*}
&b_j(x_j)=\frac{1}{2}\mu_2(\mathcal{K})(\beta_j^*)''(x_j)\rho_{Kj2}+o_p\left(\rho_{Kj2}\right),\\
&v_j(x_j)=\frac{R(\mathcal{K})\sigma_j^2(x_j)\rho_{Kj,-1}}{N_Kp_j(x_j)}+o_p\left(N_K^{-1}\rho_{Kj,-1}\right),
\end{align*}
where $\sigma_j^2(x_j)$ is defined as in \eqref{eq:sigma} and $m^*(\bm{x})=g^{-1}(\beta^*(\bm{x}))$.
\end{thm}

The proposed online estimate is asymptotically oracle in the sense that the asymptotic distribution of $\widetilde\beta_{Kj}(x_j)$ is the same as the online estimate of $\beta_j^*(x_j)$ when $\sum_{l\neq j}\beta_l^*(x_l)$ is known.
One can derive 
\begin{align*}
IMSE(\widetilde{\beta}_{Kj})&=\int \{b_j^2(x_j)+v_j(x_j)\} dx_j\\
&=\frac{1}{4}\theta_j\rho_{Kj,2}+\frac{R(\mathcal{K})\sigma_j^2\rho_{Kj,-1}}{N_K}+o_p\left(\rho_{Kj,2}+N_K^{-1}\rho_{Kj,-1}\right),
\end{align*}
where $\theta_j,\sigma_j^2$ are defined in \eqref{eq:theta sigma}. Let $\widehat{\bm{\beta}}_K$ be the batch competitor which satisfies
\begin{eqnarray*}
IMSE(\widehat{\beta}_{Kj})=\frac{1}{4}\theta_j\widehat{h}_{Kj}^4+\frac{R(\mathcal{K})\sigma_j^2}{N_K\widehat{h}_{Kj}}+o_p\left(\widehat{h}_{Kj}^4+N_K^{-1}\widehat{h}_{Kj}^{-1}\right),
\end{eqnarray*}
where $\widehat{h}_{Kj}$ is the estimated bandwidth by batch method.
We introduce the relative efficiency
\begin{equation*}
eff(\widetilde{\beta}_{Kj})=IMSE(\widehat{\beta}_{Kj})/IMSE(\widetilde{\beta}_{Kj})
\end{equation*}
to measure performance of the proposed online method compared to the classical estimate using full data.
There would be no efficiency loss if $\widetilde{\bm{\eta}}_{k|K}=\widehat{\bm{h}}_{K}$ for all $k=1,2,\ldots,K$. Hence the intuition is to make $\widetilde{\bm{\eta}}_{k|K}$ as close to $\widehat{\bm{h}}_{K}$ as possible for all $k$. 
The difference between $\widetilde{\bm{\eta}}_{k|K}$ and $\widehat{\bm{h}}_{K}$ is determined by two factors: the selection of $\widetilde{h}_{Kj}$ and the form of $\{\eta_{kj\ell}\}_{\ell=1}^L$. 
In parallel with ${(\widehat{h}_{Kj}-h_{Kj}^*)}/{h_{Kj}^*}=O_p(N_K^{-2/7})$  in \cite{mammen_bandwidth_2005},
we present the convergence of the online bandwidth $\widetilde{\bm{h}}_{K}$ as follows.

\begin{thm}
\label{THM:BAND}
  Let $\widetilde{h}_{Kj}$ be the proposed online bandwidth \eqref{eq:online h}. If the bandwidths $h_{Kj}^\theta,h_{Kj}^\sigma$ and candidates $\eta_{Kjl}^{\theta},\eta_{Kjl}^{\sigma}$ for pilot estimates satisfy \eqref{eq:gr} and \eqref{eq:candidate r} for $j=1,2,\ldots,d$, 
  then as $K\rightarrow\infty$, under Assumptions \eqref{assump: density}-\eqref{assump: pseudo-bandwidths},
  \begin{equation*}
  \frac{\widetilde{h}_{Kj}-h_{Kj}^*}{h_{Kj}^*}=O_p\left(N_K^{-\frac{2}{7}}\right),
  \end{equation*}
  where $h_{Kj}^*$ is the optimal bandwidth in \eqref{eq:opt h}.
\end{thm}

Since $\widetilde{h}_{Kj}$ is strictly decreasing with respect to $K$, the candidate sequence $\{\eta_{kj\ell}\}_{\ell=1}^L$ shall be no larger than $\widetilde{h}_{kj}$ so that we can backtrack an appropriate pseudo-bandwidth in future estimates. Without loss of generality, we set $\eta_{kj1}>\eta_{kj2}>\ldots>\eta_{kjL}$, and derive the form of $\eta_{kj\ell}$ by optimizing the relative efficiency  $eff(\widetilde{\beta}_{Kj})$. A lower bound of $eff(\widetilde{\beta}_{Kj})$ driven by $L$ is also obtained in the theorem below.

\begin{thm}
\label{THM:EFF} 
Under Assumption \eqref{assump: density}-\eqref{assump: pseudo-bandwidths}, when the candidate bandwidth sequence takes the form
\begin{equation}\label{eq:eta}
\eta_{kj\ell}=\left(\frac{L-l+1}{L}\right)^{1/5}\widetilde h_{kj},
\end{equation}
where $\widetilde h_{kj}$ as in \eqref{eq:online h} is the onilne estimate of optimal bandwidth whose pilot estimates satisfying \eqref{eq:gr} and \eqref{eq:candidate r}, then as $K\rightarrow\infty$, the corresponding online estimate $\widetilde\beta_{Kj}$ has the asymptotic distribution as stated in Theorem \ref{THM:AN} and attains optimal relative efficiency with the following lower bound,
\begin{equation*}
\left(1+c_1\frac{1}{L}+c_2\frac{1}{L^2}\right)^{-1}+O_p\left(\frac{N_K^{-\frac{1}{5}}}{L}+N_K^{-\frac{2}{7}}\right),
\end{equation*}
where $c_1=0.183$ and $c_2=0.003$.
\end{thm}

\begin{figure}[htbp]
  \centering
  \includegraphics[width = 3.5in]{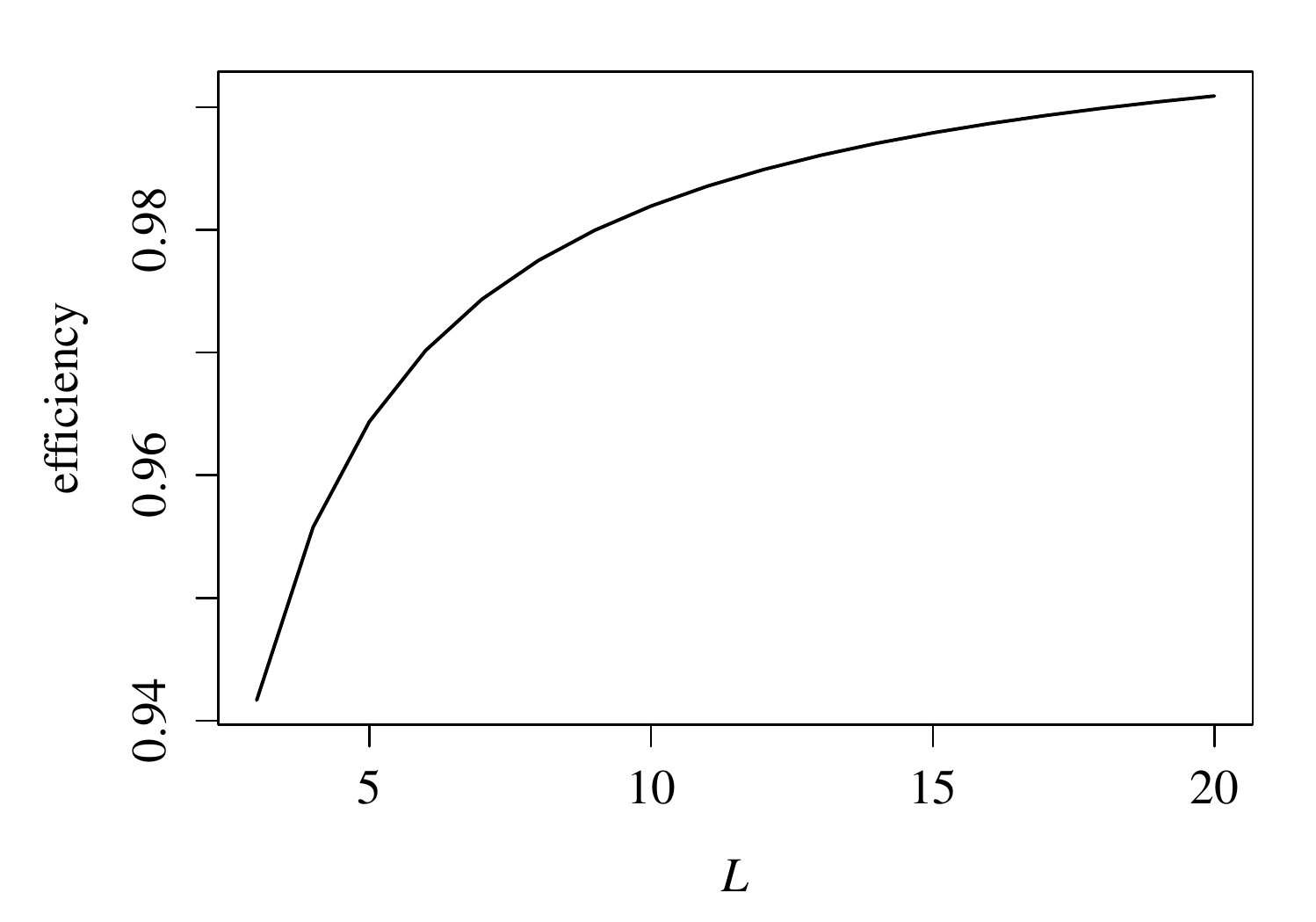}
  \caption{\label{fig:theoretical_eff} The theoretical lower bound for the relative efficiency of the proposed online estimates versus different lengths $L$ of candidate bandwidth sequences.}
\end{figure}

This bound is the same as the one-dimensional nonparametric regression discussed in \cite{yang2021ofda}. From the illustration of Figure \ref{fig:theoretical_eff}, the relative efficiency of the proposed estimate improves rapidly as $L$ increases, and exceeds 95\% when $L\ge5$. This is desirable to attain high efficiency with a small/moderate $L$. Note that the computational cost in terms of time and memory is proportional to $L$, this lower bound helps make an informed
trade-off between statistical and computational efficiency, which makes the proposed method
practically useful.
{}
In Theorem \ref{THM:AN}--\ref{THM:EFF}, we discussed the properties of $\widetilde{\bm{\beta}}_K$ Now we delineate the convergence of the algorithm, i.e., $\widetilde{\bm{\beta}}_K^{(m)}\rightarrow\widetilde{\bm{\beta}}_K$ as $m\rightarrow\infty$. Let $B_r(\widetilde{\bm{\beta}}_K)$ denote the ball with the center $\widetilde{\bm{\beta}}_K$ and a radius $r$. We state the following results for the  outer and inner iterations of the proposed algorithm.
\begin{thm}
\label{THM:OUTER CONVERGENCE}
Let $\widetilde{\bm{\beta}}_K^{(m)}$ be the $m$th outer step estimate of \eqref{eq:linearapprox} at time $K$. Under assumptions \eqref{assump: density}-\eqref{assump: pseudo-bandwidths}, there exist fixed $r$, $C>0$ and $0<\gamma<1$ that satisfy: if the initial values $\widetilde{\bm{\beta}}_K^{(0)}\in B_r(\widetilde{\bm{\beta}}_K)$ with probability tending to 1, then
\begin{equation*}
P\left(\Vert\widetilde{\bm{\beta}}_K^{(m)}-\widetilde{\bm{\beta}}_K\Vert_p\leq C2^{-(m-1)}\gamma^{2^k-1}\right)\rightarrow1.
\end{equation*}
\end{thm}

\begin{thm}
\label{THM:INNER CONVERGENCE}
Under assumptions \eqref{assump: density}-\eqref{assump: pseudo-bandwidths}, the inner iteration \eqref{inner iter} converges at a geometric rate. If the initial values $\widetilde{\bm{\beta}}_K^{(0)}\in B_r(\widetilde{\bm{\beta}}_K)$ with probability tending to 1, then the geometric convergence of the inner iteration is uniform for all steps in the outer iteration, with probability tending to one.
\end{thm}

For the initial values of Algorithm \ref{alg:ogam}, when $K=1$, we set $\bm{\beta}^{(0)}=\widetilde{\bm{\beta}}_0$ where $\widetilde{\bm{\beta}}_0$ is the parametric model fit following \cite{yu_smooth_2008}, and when $K>1$, we use the natural initialization $\bm{\beta}^{(0)}=\widetilde{\bm{\beta}}_{K-1}$. This may not be contained in the ball $B_r(\widetilde{\bm{\beta}}_K)$ with probability tending to one. As suggested in \cite{yu_smooth_2008}, an alternative is to use the marginal integration estimator proposed by \cite{linton1996estimation}, which is consistent but computationally expensive. As our focus is the comparison between the online and batch methods, we follow \cite{yu_smooth_2008} and use the parametric model fit as initials.

\section{Numerical Experiments}
\label{sec:numerical}
\subsection{Simulation}

We conduct simulation to illustrate the performance of the proposed online method and
verify the theoretical findings in Section \ref{sec:th}.
The simulation is done under the following model for the conditional distribution: $Y\mid\bm{X}=\bm{x}\sim\textrm{Poisson}(m(\bm{x}))$, where
\begin{equation*}
\log(m(\bm{x}))=2+\beta_1(x_1)+\beta_2(x_2)=2+\cos(\sqrt{2}\pi x_1)+\{x_2+\sin(2\pi x_2)\}/2.
\end{equation*}
The covariate vector
$\bm{X}=(X_1,X_2)\sim N_2(0,0;10,10;0.9)\textrm{ truncated on }[0,1]^2$,
where $N_2(\mu_1,\mu_2; \sigma_1^2,\sigma_2^2,\rho)$ denotes the bivariate normal distribution with means $\mu_1,\mu_2$, variances $\sigma_1^2,\sigma_2^2$, and correlation coefficient $\rho$.
The components satisfy the normalizing constraint given at \eqref{eq:constrain} are $\beta_1^*(x_1)=\cos(\sqrt{2}\pi x_1)-0.227$, $\beta_2^*(x_2)=\{x_2+\sin(2\pi x_2)\}/2-0.317$ and $\beta_0^*=2.544$.

The sample size of each block $n_k$ (before rounding to the nearest integer) is normally distributed with mean 100 and standard deviation 10.
We have experimented an extensive range of $G$ and $R$ in \eqref{eq:gr} for pilot estimates of $\theta_j$ and $\sigma_j^2$ and find that the convergence rate of the online bandwidth is not sensitive to the values of $G$ and $R$. Setting them in $[0.5,1]$ is in general adequate, thus we set $G=R=0.5$ in the sequel. 
Also note that the relative efficiency increases rapidly to 1 from Figure \ref{fig:theoretical_eff}, it suffices to set $L'=10$ as the length of candidate bandwidth sequences for the pilot estimates of $\theta_j$ and $\sigma_j^2$, $j=1,2$. The experiment is replicated 100 times, each with 1000 blocks. 

\begin{figure}[htbp]
  \centering
  \includegraphics[width = 5.2in]{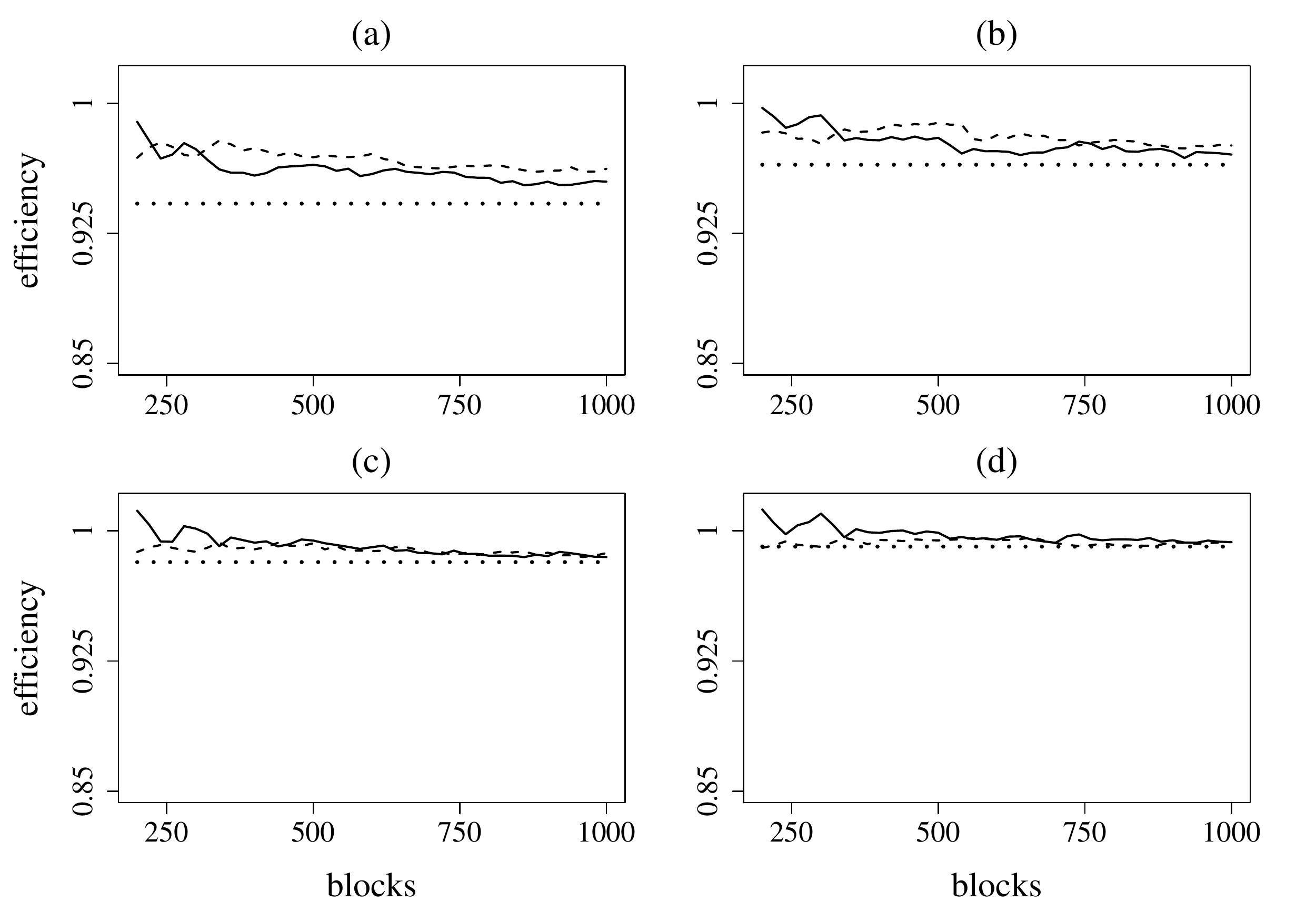}
  \caption{\label{fig:empirical_eff} The empirical relative efficiency for components $\beta_1$ (solid) and $\beta_2$ (dashed) the theoretical lower bound (dotted) of the proposed estimate using $L=3,5,10$ and $20$,  corresponding to (a)-(d), respectively.}
\end{figure}

\begin{figure}[htbp]
  \centering
  \includegraphics[width = 5.2in]{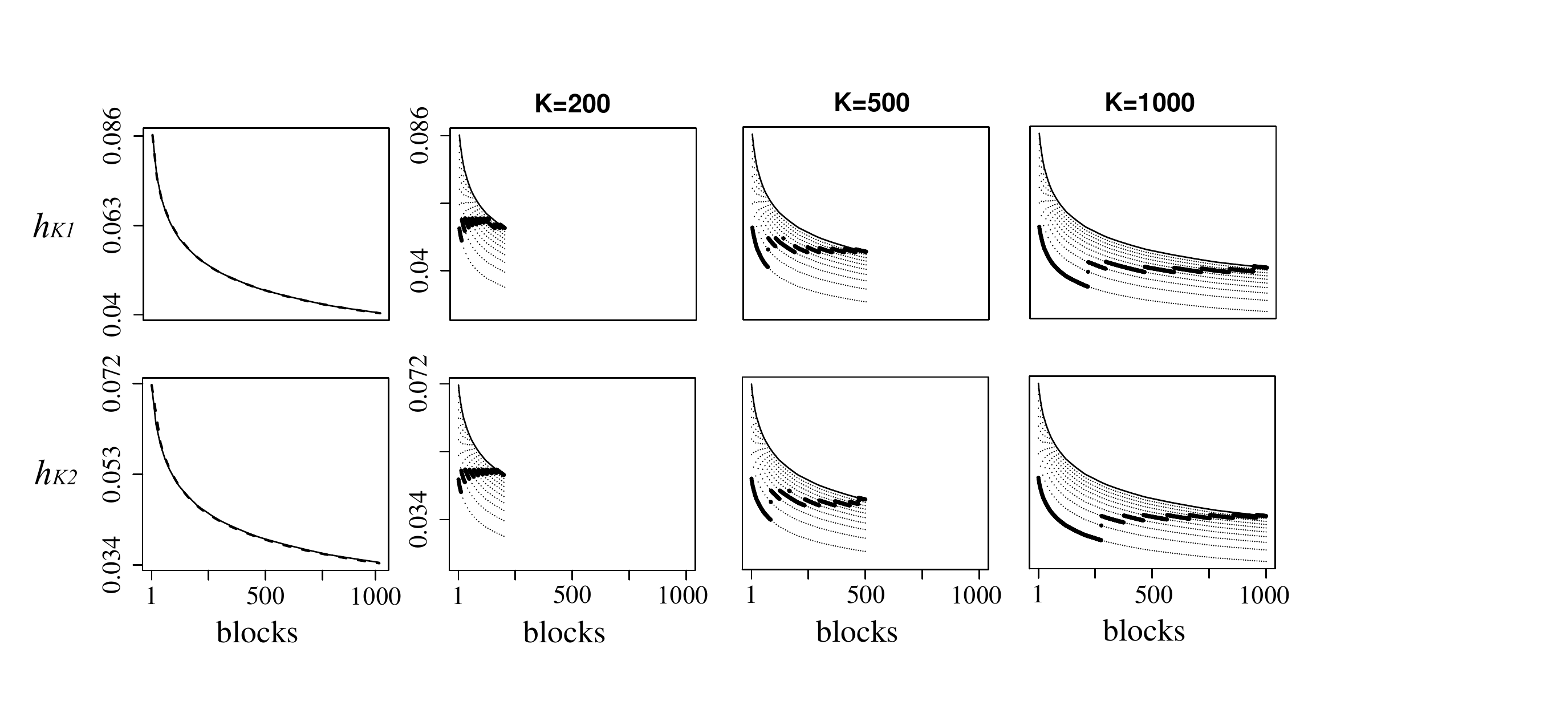}
  \caption{\label{fig:band} The two rows show, respectively, the bandwidth selection for estimating $\beta_1$ and $\beta_2$. The first column shows the Monto Carlo averages of the bandwidths selected by the online method (solid) $\widetilde{h}_{kj}$ and the batch method $\widehat{h}_{kj}$ (dashed) based on 100 runs for $k=1,2,\ldots,K$, and the other three columns depict the dynamically updated bandwidths $\{\widetilde{\eta}_{k|K,j}\}_{k=1}^K$ (thick dots) at time $K=200,500,1000$, respectively, along with the dynamic candidate sequences (light dots) and optimal bandwidth estimates (connected by solid line) at each $k$.}
\end{figure}

\begin{figure}[htbp]
  \centering
  \includegraphics[width = 3.5in]{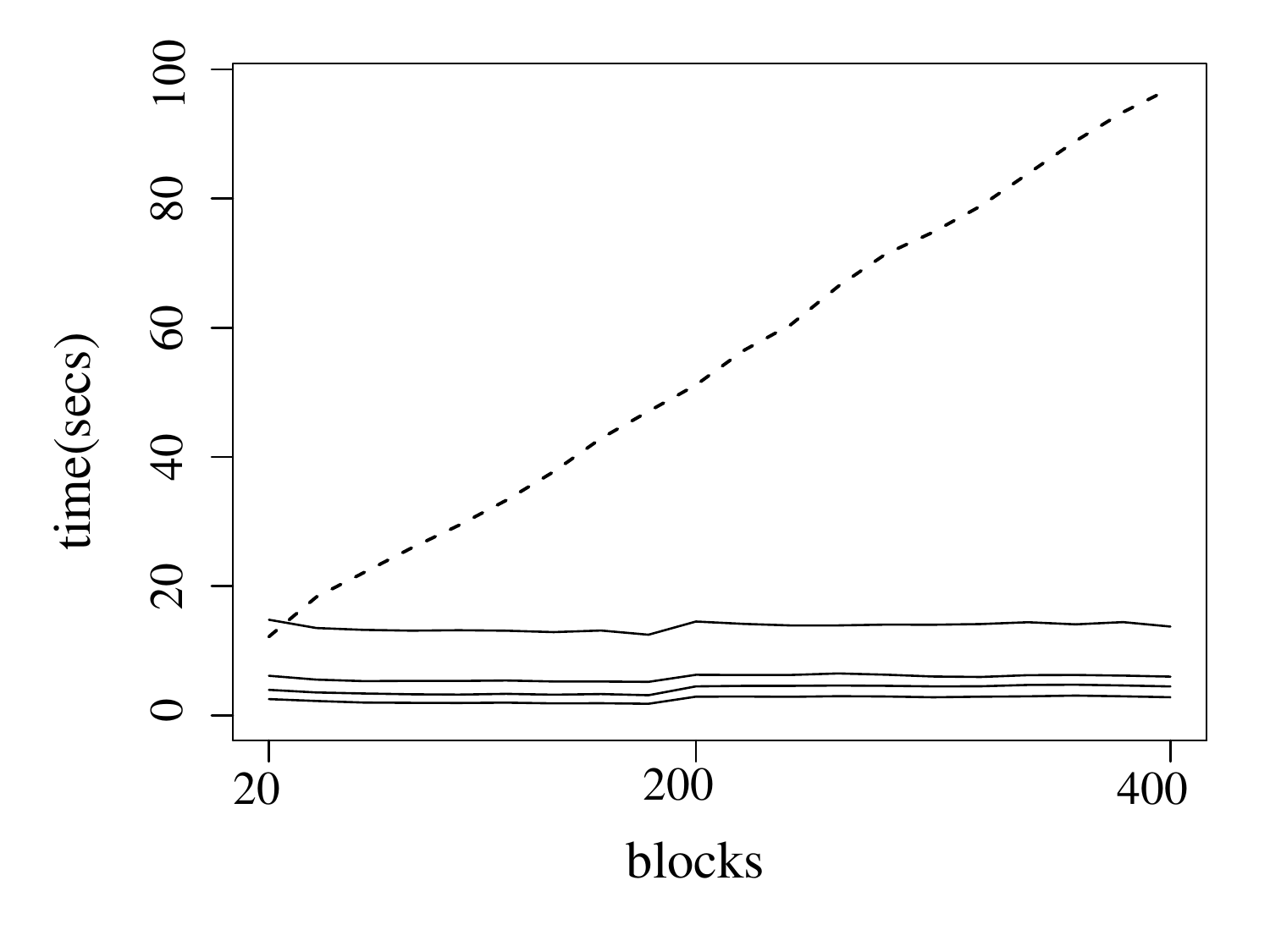}
  \caption{\label{fig:simutime} The comparison of computing times between our method (solid) and the batch method (dashed), and the four dashed lines from the lowest to the highest in the left panel represent $L=3,5,10, 20$, respectively.}
\end{figure}

Figure \ref{fig:empirical_eff} shows the empirical relative efficiency that increases with $L$ and is stably higher than the theoretical lower bound in Theorem \ref{THM:EFF} when $K$ tends large. 
The convergence of bandwidth in Theorem \ref{THM:BAND} is also examined in Figure \ref{fig:band}, where both batch and online selections converge to the theoretical optimal bandwidth as data accumulate. Moreover, we also depict the dynamically updated bandwidths $\{\widetilde{\bm{\eta}}_{k|K}\}_{k=1}^K$ that are  used to produce the estimates at time $K=200,400,1000$, similar to that shown in Figure \ref{fig:sketch}. This provides empirical support for the fact that the proposed dynamic updating algorithm is indeed able to adjust the sub-sufficient statistics in previous blocks in spite of no access to those data. Lastly, we compare the computing times of the batch and the proposed online methods on a Unix server of 2.10GHz CPU and 188G memory with 176 logic cores. It is noted from Figure \ref{fig:simutime} that the computing time of the batch method increases approximately linearly as data accumulate, and the online method using $L=3,5,10, 20$ spends nearly constant times proportional to $L$. Thus we graphed the first 400 blocks for better visualization and remark that the online method achieves substantial computational saving. We close the simulation study by suggesting that, based on our empirical and theoretical findings, a reasonable range for choosing $L$ is $[5,20]$, depending on whether accuracy or computation is of main concern.

\subsection{Real data application}

In this subsection, we present the airline delay example to illustrate the usefulness of the proposed online method. The dataset consists of flight arrival and departure details for all commercial airports in the USA (\url{https://community.amstat.org/jointscsg-section/dataexpo/dataexpo2009}). We use the data of flights departing during 6:00 to 23:00 from January 1996 to December 2004 to fit a logistic regression and model the probability of late arrival $\pi$. We select scheduled departure time $X_1$ and historical delay rate of each flight $X_2$ as the covariates. We mention that in the raw dataset, the variable $\tt{FlightNum}$ does not contain the information of routines of flights, which results in duplicate flight numbers for each day. We combine the original variables $\tt{Origin}$ and $\tt{FlightNum}$ as the new flight number to compute the number of historical travels and historical delay rate. Then we select the flights which have more than 30 historical travels and whose historical delay rates lie between 15\% and 65\% to train the model. The observations are divided into blocks by date, i.e., set data of a day as a block. Finally there are 3283 blocks in total, and each block contains 3339--8776 observations.
The parameters are set as $L=L'=10$ and $G=R=1$.
To mitigate the heavy tail of $X_2$, we make the following transformation of $X_2$: $X_2'=\log(X_2-0.14)$.
Then both $X_1$ and $X_2'$ are normalized to be on $[0,1]$. 

\begin{figure}[htbp]
  \centering
  \includegraphics[width = 5in]{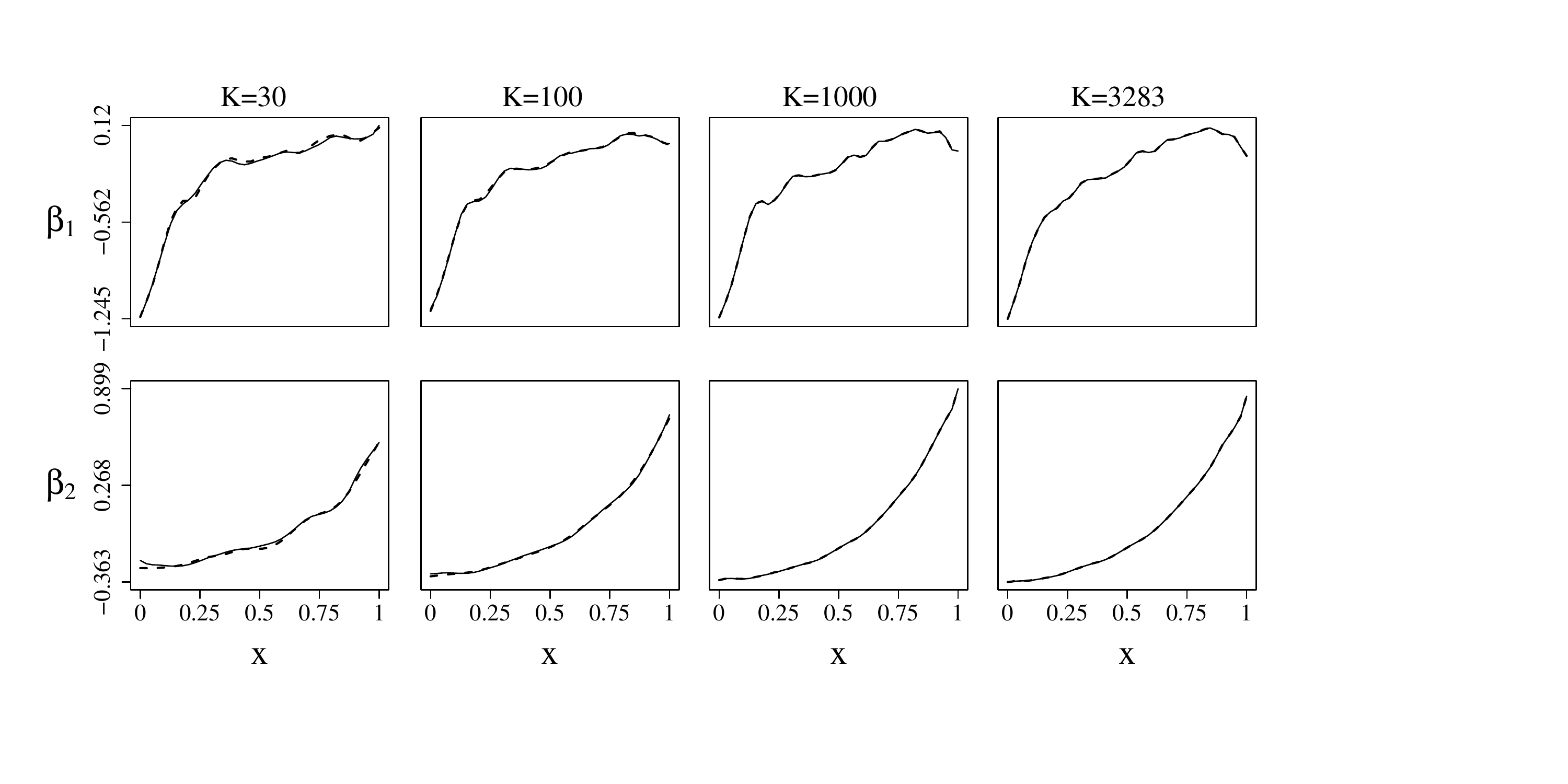}
  \caption{\label{fig:flights betas} The estimates of $\beta_1$ and $\beta_2$ using the proposed online method in solid lines ($\widetilde\beta_{K1},\widetilde\beta_{K2}$) and the classical batch method in the dashed lines ($\widehat\beta_{K1},\widehat\beta_{K2}$) at time $K=30,100,1000$ and 3283.}
\end{figure}

\begin{figure}[htbp]
  \centering
  \includegraphics[width = 3.5in]{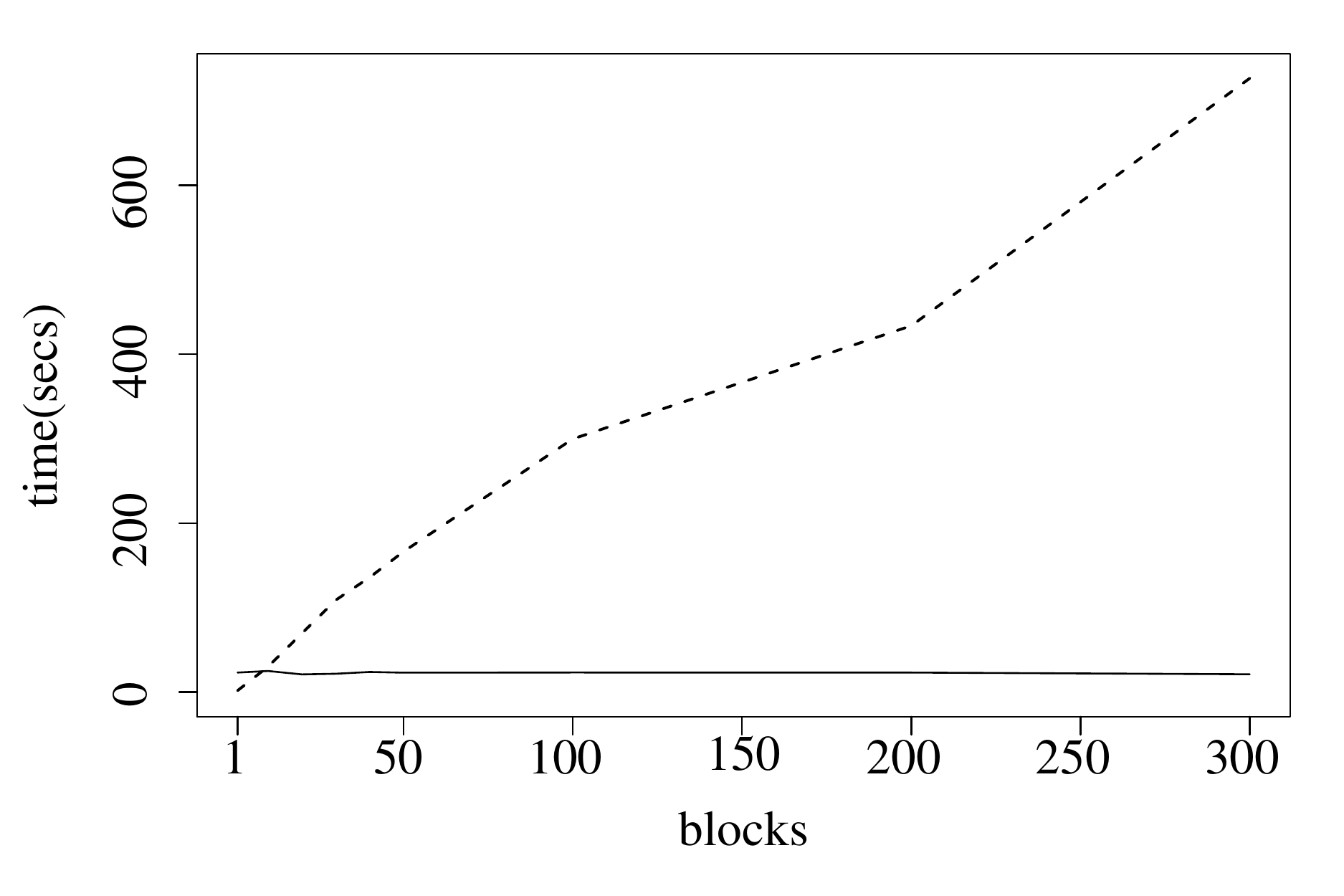}
  \caption{\label{fig:flights time} The comparison of computing times between our method (solid) and the batch method (dashed) on the flights data set.}
\end{figure}

Denote the component functions as $\beta_1(X_1)$ and $\beta_2(X_2)$, i.e., $\textrm{logit}(\pi)=\beta_0+\beta_1(X_1)+\beta_2(X_2)$, where $\textrm{logit}(\pi)=\log(\pi(1-\pi)^{-1})$. Figure \ref{fig:flights betas} present the proposed online estimates $\widetilde\beta_{K1}(\cdot),\widetilde\beta_{K2}(\cdot)$ and the batch competitors $\widehat\beta_{K1}(\cdot),\widehat\beta_{K2}(\cdot)$ when $30,100,1000,3283$ blocks are included.
The probability of delay increases along scheduled departure time and historical delayed rate.
It can be seen that $\widetilde\beta_{Kj}$ converges to $\widehat\beta_{Kj}$ for $j=1,2$ as data accumulate, which validates the usefulness of our method. Then we use the data of year 2005 as the test set which contains $N_{test}=2998231$ flights. The $i$-th flight is predicted to be delayed $\widetilde Y_i=1$ (or $\widehat Y_i=1$) if $\widetilde\pi_i>0.5$ (or $\widehat\pi_i>0.5$).  The prediction errors $\widetilde{err}=\#\{\widetilde Y_i\neq Y_i\}/N_{test}$ and $\widetilde{err}=\#\{\widehat Y_i\neq Y_i\}/N_{test}$ at different time points are listed in Table \ref{tab:pre error}. The prediction errors using online and batch methods both decrease to 0.238 as more data enter the model. Finally, we compare the computing times of the online and batch methods in Figure \ref{fig:flights time} till $K=300$ which again verifies that our online method is more computationally efficient.

{\small\begin{table}[htbp]
\centering
\caption{Prediction Errors of Flights Dataset}
  \begin{tabular}{|c|r|r|r|r|r|r|r|r|r|}
  \hline
  $K$     & \multicolumn{1}{c|}{1} & \multicolumn{1}{c|}{10} & \multicolumn{1}{c|}{30} & \multicolumn{1}{c|}{50} & \multicolumn{1}{c|}{100} & \multicolumn{1}{c|}{500} & \multicolumn{1}{c|}{1000} & \multicolumn{1}{c|}{2000} & \multicolumn{1}{c|}{3283} \\
  \hline
  online & 0.529 & 0.361 & 0.239 & 0.239 & 0.239 & 0.238 & 0.238 & 0.238 & 0.238 \\
  \hline
  batch & 0.529 & 0.385 & 0.239 & 0.239 & 0.238 & 0.238 & 0.238 & 0.238 & 0.238 \\
  \hline
  \end{tabular}%
\label{tab:pre error}%
\end{table}%
}

\section*{Supplementary Materials}

The supplementary material contains the proofs for Theorem \ref{THM:AN}--\ref{THM:INNER CONVERGENCE}.
\par
\section*{Acknowledgements}

This research is supported by National Natural Science
Foundation of China Grants No.11931001 and 11871080, the LMAM, and the Key Laboratory of Mathematical Economics and Quantitative Finance (Peking University), Ministry of Education.
\par


\bibhang=1.7pc
\bibsep=2pt
\fontsize{9}{14pt plus.8pt minus .6pt}\selectfont
\renewcommand\bibname{\large \bf References}
\expandafter\ifx\csname
natexlab\endcsname\relax\def\natexlab#1{#1}\fi
\expandafter\ifx\csname url\endcsname\relax
\def\url#1{\texttt{#1}}\fi
\expandafter\ifx\csname urlprefix\endcsname\relax\def\urlprefix{URL}\fi

\bibliographystyle{chicago}      
\bibliography{reference}   



\end{document}